\documentclass[aop]{imsart}

\RequirePackage{amsthm,amsmath,amsfonts,amssymb}
\RequirePackage[numbers]{natbib}
\RequirePackage[colorlinks,citecolor=blue,urlcolor=blue]{hyperref}
\RequirePackage{graphicx}

\startlocaldefs
\theoremstyle{plain}

\newtheorem{theorem}{Theorem}[section]
\newtheorem{lemma}[theorem]{Lemma}
\theoremstyle{remark}

\newcommand{\dif}{\mathrm{d}}
\endlocaldefs

\begin{document}

\begin{frontmatter}
\title{The Gaussian free-field as a stream function:\\
asymptotics of effective diffusivity in infra-red cut-off}
\runtitle{Gaussian free-field as stream function}

\begin{aug}
\author[A]{\fnms{Georgiana}~\snm{Chatzigeorgiou}\ead[label=e1]{georgiana.chatzigeorgiou@mis.mpg.de}},
\author[A]{\fnms{Peter}~\snm{Morfe}\ead[label=e2]{morfe@mis.mpg.de}},
\author[A]{\fnms{Felix}~\snm{Otto}\ead[label=e3]{felix.otto@mis.mpg.de}},
\and
\author[A]{\fnms{Lihan}~\snm{Wang}\ead[label=e4]{lihan.wang@mis.mpg.de}}
\address[A]{Max Planck Institute for Mathematics in the Sciences \printead[presep={,\ }]{e1,e2,e3,e4}}

\end{aug}

\begin{abstract}
We analyze the large-time asymptotics of a passive tracer with drift equal to the curl of the Gaussian free field in two dimensions with ultra-violet cut-off at scale unity.  We prove that the mean-squared displacement scales like $t \sqrt{\ln t}$, as predicted in the physics literature.  This improves the asymptotics recently established in recent work of Cannizzaro, Haunschmidt-Sibitz, and Toninelli (2022), which uses mathematical-physics type analysis in Fock space.  Our approach involves studying the effective diffusivity $\lambda_{L}$ of the process 
with an infra-red cut-off at scale $L$, and is based on techniques from stochastic homogenization.
\end{abstract}



\end{frontmatter}

\section{Setting and statement of main result}\label{s:intro}

The aim of this paper is to analyze, via homogenization methods, the large-scale behavior of a passive tracer subject to diffusion and convection by a divergence-free drift $b$ in $\mathbb{R}^d$. 
More precisely, we assume that $b$ is a stationary centered Gaussian vector field and as such characterized by its covariance function $R(x - y) = \mathbb{E}b(x) \otimes b(y)$, or rather the Fourier transform thereof\footnote{Here $\iota = \sqrt{-1}$ is the imaginary unit.} ${\mathcal{F}R}(k) = (2\pi)^{-1}\int \dif x \exp(-\iota k \cdot x) R(x)$, by imposing
\begin{align} \label{E: spectrum of b}
    {\mathcal{F}}R(k) = \text{Id} - |k|^{-2} k \otimes k.
\end{align}
This is the most natural ensemble for a divergence-free drift, since it is the Leray projection of white noise. White noise is the only Gaussian process with a finite range of dependence, and from an information-theoretical point of view has the most randomness.

\medskip

The passive tracer is modelled via the stochastic differential equation (SDE)
\begin{align*}
    \dif X_{t} =  b(X_{t}) \, \dif t + \sqrt{2} \dif W_{t}.
\end{align*}
Here $W_t$ is the standard $d$-dimensional Brownian motion. It is well-known that Brownian motion $W_t$ is invariant in law under the parabolic rescaling $\ell^{-1} W_{\ell^2 t}$.
Since velocities scale like $\ell^{-1}$ under parabolic rescaling, the solution $X_t$ of the SDE has the same invariance if and only if the drift $b$ is invariant in law under $\ell b(\ell\cdot)$. Since our white-noise type ensemble is invariant in law under $\ell^\frac{d}{2} b(\ell\cdot)$ we see that this is the case if and only if $d=2$.  Hence this is the critical  and interesting dimension. We therefore restrict ourselves to the $d=2$ setting.

\medskip

However, the aforementioned SDE is not well-posed due to the roughness of $b$. In fact, it is critical in a stricter sense of singular SPDEs, 
namely in terms of the forward equation, which is the heat equation with a random drift (see \eqref{E:main_operator}); as it turns out, the latter is not renormalizable. Hence it is customary to introduce an ultra-violet (UV) cut-off on $b$. We may non-dimensionalize such that the UV cut-off scale is unity, so that under such modification, the covariance of $b$ now has Fourier transform
\begin{align} \label{e:no-infrared}
{\mathcal{F}}R(k) = I(|k|\le 1)(\text{Id} - |k|^{-2} k \otimes k).
\end{align}
This problem setting comes with a non-dimensional\footnote{The parameter $\epsilon$ is dimensionless since, in view of \eqref{e:no-infrared}, $b$ itself is dimensionless.} parameter,
which we implement as a prefactor $\epsilon$ of the drift;
$\epsilon$ can be assimilated to the P\'eclet number\footnote{Since the UV cut-off at scale unity provides a characteristic length scale for the drift, the P\'{e}clet number is $\epsilon \left( \int_{\mathbb{R}^{2}} \dif k |k|^{-2} \text{tr}(\mathcal{F}R(k)) \right)^{\frac{1}{2}}$.}  in the jargon of transport-diffusion processes:
	\begin{equation} \label{e:sde}
		\dif X_{t} = \epsilon b(X_{t}) \, \dif t + \sqrt{2} \dif W_{t}.
	\end{equation}
   While it is customary to fix $\epsilon^{2} = 1$,\footnote{Since $b$ and $-b$ have the same law, the annealed law of $X_{t}$ depends only on $\epsilon^{2}$.} the presence of $\epsilon^{2}$ not only provides a convenient ordering principle but we also use a monotonicity in $\epsilon^{2}$ (cf.\ Lemma \ref{l:monotone}).

\medskip

Our main result amounts to capturing the borderline super-diffusive behavior $\mathbb{E}|X_{t}|^{2} \sim t \sqrt{\ln t}$ of the solution to the SDE \eqref{e:sde}. Due to statistical isotropy of the problem setting, $\mathbb{E}(\xi\cdot X_{t})^2$ does not depend on the choice of unit vector $\xi$.
\begin{theorem}\label{t:laplace}
	For every $\epsilon^2 < \infty$, and every unit vector $\xi$, there exists a constant $C_{\epsilon}$ such that for every $t \geq 1$ 
	\begin{equation}\label{eq:EXT2}
C_{\epsilon}^{-1}	\le	\frac{\mathbb{E}(\xi\cdot X_{t})^2}{2t \sqrt{1+\frac{1}{2}\epsilon^2 \ln t}}  \le C_{\epsilon}.
	\end{equation}
Moreover, $C_{\epsilon} =1+ O(|\epsilon|)$ as $\epsilon^2\to 0$.
\end{theorem}
For our analysis we introduce an additional infra-red (IR) cut-off at a scale $L$ and study the asymptotics as $L \to \infty$. More precisely, let $b_L$ denote the centered stationary Gaussian vector field with covariance function determined by its Fourier transform through 
\begin{align} \label{e:infrared}
{\mathcal{F}}R_L(k) = I(L^{-1} \le |k|\le 1)(\text{Id} - |k|^{-2} k \otimes k).
\end{align}
In this setting, the effective behavior of 
\begin{equation*} 
\dif X_t^{(L)} = \epsilon b_L(X_t^{(L)}) \dif t + \sqrt{2}\dif W_t,
\end{equation*}
is known to be diffusive (see Section \ref{s:heur} for more details and \cite{KLO12} for a general theory), in the sense that the following limit, known as the effective diffusivity
\begin{equation} \label{e:diffusivity-X_L}
   \lambda_L:= \lim_{t\to \infty} \frac{\mathbb{E}(\xi\cdot X_{t}^{(L)})^2}{2t} 
\end{equation}
exists. The core of our argument concerns the scaling of $\lambda_L$ in terms of $L\gg 1$. To see why the scaling of $\lambda_L$ is relevant to the scaling of $\mathbb{E}|X_t|^2$, we draw inspiration from the heuristics
suggested by Fannjiang in \cite{F98}: as long as $\lambda_L t \ll  L^2$, 
the environment on length scales larger than $L$ has little effect on
the process $X_t$, 
hence at least heuristically, the processes $X_t$ and $X_t^{(L)}$ should have similar behavior. The combination of Theorems \ref{t:laplace} and \ref{t:main} give a rigorous justification of this heuristics. 

\medskip

Our results will be obtained by deriving the large-scale, i.e., ``effective" or homogenized behavior, of the generator of the process  \eqref{e:sde}
\begin{equation}\label{E:main_operator}
		-\Delta u - \epsilon b \cdot \nabla u.
	\end{equation}
 Combining \eqref{e:no-infrared} and \eqref{e:infrared}, we observe that \begin{equation}\label{e:bminbL} \begin{aligned}
    \mathbb{E}|b - b_L|^{2} =\mathrm{tr}(R-R_L)(0) & = \int \frac{\dif k}{2\pi} \, \mathrm{tr}\mathcal{F}(R-R_L)(k)  \\ & = \int_{0\le |k| \le L^{-1}} \frac{\dif k}{2\pi} \, \mathrm{tr}\Big(\mathrm{Id}- \frac{k\otimes k}{|k|^2}\Big) = \frac{1}{2L^2},
\end{aligned}\end{equation}  so that $b_L$ is at least heuristically an  approximation of $b$ as $L\to\infty$. Therefore, our strategy consists in approximating \eqref{E:main_operator} via the associated operator
\begin{align}\label{ao36}
-\Delta u-\epsilon b_L \cdot \nabla u.
\end{align}

\medskip

Analogously to \eqref{e:diffusivity-X_L}, the effective behavior of \eqref{ao36} is described by the isotropic diffusion operator $-\lambda_L\Delta$. This is an instance of convection-enhanced diffusion, since qualitatively $\lambda_L>1$
for $L>1$ (while trivially $\lambda_{L=1}=1$), see (\ref{ao08}). The presence of a widening spectrum of scales explores the limits of homogenization theory, which typically involves a scale separation. For $\epsilon^2=1$ and in the regime of $L\gg 1$, it is conjectured that $\lambda_L^2\sim\ln L$, see \cite{F98}. We establish a more precise version of this:

\begin{theorem}\label{t:main}
   For every $\epsilon^2 < \infty$, there exists a constant 
$C_\epsilon$ such that for all $L\ge 1$
    \begin{equation}\label{ao38}
        C_\epsilon^{-1} \le \frac{\lambda_L^2}{1+\epsilon^2\ln L} \le C_\epsilon.
    \end{equation}
Moreover, $C_{\epsilon} =1+ O(|\epsilon|)$ as $\epsilon^2\to 0$.
\end{theorem}

Note that while our result correctly captures the asymptotics in the regime $\epsilon^2\ll 1$, its proof provides no information 
on the interesting asymptotics in the opposite regime $\epsilon^2\gg 1$.

\medskip

A key fact which we use extensively in what follows is that in $d=2$, the drift $b_L$ can be expressed as the curl of the so-called stream function $\psi_L$, or equivalently $b_L= -J\nabla \psi_L$, where $J$ denotes the rotation by 90 degrees
\begin{align}\label{ao95}
J=\left(\begin{array}{cc}0&-1\\1&0\end{array}\right)\quad\mbox{so that
by symmetry of second derivatives}\quad\nabla\cdot J\nabla=0.
\end{align}
As a result, the drift-diffusion operator (\ref{ao36}) can be rewritten in divergence form:
\begin{align}\label{ao25}
-\Delta u-\epsilon b_L \cdot \nabla u=-\nabla\cdot a_L\nabla u
\quad\mbox{where}\quad a_L:=1+\epsilon\psi_LJ.
\end{align}

Due to the form of $b_L$, the stream function $\psi_L$ is the two-dimensional Gaussian free field (GFF) with UV cut-off at scale 1 and IR cut-off at scale $L$. More precisely, $\psi_L$ is the centered stationary Gaussian scalar field with covariance function $c_L$ determined through
\begin{align}\label{ao39}
{\mathcal F}c_L(k)=I(L^{-1}\le|k|\le 1)|k|^{-2}.
\end{align}
Since $\mathcal{F}c_{L}(k)$ only depends on the length $|k|$ of the wave vector $k$,
$c_{L}(x-y)$ only depends on $|x-y|$, which implies that the centered Gaussian field $\psi_{L}$
is statistically isotropic, which means that for any orthogonal $R$,
$\psi_{L}(R\cdot)$ and $\psi_L$ follow the same law.  The logarithmic character of the GFF in two space-dimensions becomes apparent when computing
\begin{align}\label{ao37}
\mathbb{E}\psi_L^2=c_L(0)=\int\frac{\dif k}{2\pi}{\mathcal F}c_L
\stackrel{(\ref{ao39})}{=}\int_{\{L^{-1}\le|k|\le 1\}}\frac{\dif k}{2\pi}|k|^{-2}=\ln L.
\end{align}

\medskip
	
The paper is organized as follows: In Section \ref{s:relwork} we give a brief overview of 
the related literature. In Section \ref{s:heur} we heuristically derive an ODE for
$\lambda_L$ as a function of the scale $L$, which by integration yields
$\lambda_L^2=1+\epsilon^2\ln L$,  and is based on homogenization theory 
and the notion of a corrector $\phi_L$.  In Section \ref{s:sket} we recall the notion of a flux corrector $\sigma_L$ designed
to estimate the homogenization error. Section \ref{s:constr} is at the core of our proof and introduces stationary approximations $(\tilde\phi_L,\tilde\sigma_L)$ of $(\phi_L,\sigma_L)$;  these proxies are estimated in Section \ref{s:approxcorr}.  In Section \ref{s:error} we estimate the error $f_L$ in this approximation. The proofs of Theorem \ref{t:main} and Theorem \ref{t:laplace} are deferred to Sections \ref{s:proofmain} and \ref{s:proofmain2}, respectively.



\section{Related work}\label{s:relwork}

Mathematically rigorous homogenization results for passive tracers convected by  
a divergence-free drift date back to the work of 
Papanicolaou and Varadhan \cite{PV81}, Osada \cite{osada} and Oelschl\"ager \cite{oelschlager}. 
The enhancement effect of the drift was quantified for the first time in \cite{FP94} in a periodic setting. The random setting was considered for instance in \cite{FP96}, \cite{FK97}, \cite{FK99}, where diffusivity of the stochastic process  is proved for a stationary and ergodic stream function $\psi$ which satisfies suitable integrability conditions. We refer the interested 
reader to Chapters 11 and 12 of \cite{KLO12} for an overview on the topic of diffusivity 
for Brownian particles driven by a divergence-free random drift, 
and \cite{KT17}, \cite{T18}, \cite{Fe20} for the state of the art.

\smallskip

We note that homogenization is much more subtle when ensembles of
non-divergence-free drifts are considered, in which case the uniform
distribution of the tracer is no longer invariant. A scaling argument reveals that
homogenization is to be expected for $d>2$; on this basis
a discrete renormalization group argument was set up in \cite{BK91}.
A fully rigorous argument dealing with the traps was given in \cite{SZ06} for $d>2$
and small P\'eclet number. The iterative approach employed here also
uses an iteration over geometrically growing scales, more like in \cite{BK91}
than in \cite{SZ06}.

\smallskip

Let us now focus on the two-dimensional case. The Fourier transform of the GFF is borderline not integrable at $k=0$ (as opposed to \eqref{ao37}),
which means that $\psi$ is not stationary (while the drift $J\nabla\psi$ is).
In this case, as was predicted in the physics literature \cite{FFQSS85}, \cite{KB89} using renormalization-group heuristics (which we refer the readers to the Appendix \ref{sec:appendix_heuristic}), the corresponding process is (borderline) super-diffusive, in the sense that 
\begin{equation*}
    \mathbb{E}|X_t|^2 \sim t\sqrt{\ln t}, \, \mbox{ as } \, t\to\infty.
\end{equation*} 
Later, in \cite{TV12} a different heuristic explanation based on a scaling argument is given, and using the Fock space approach developed in \cite{LQSY04}, the authors rigorously establish the (sub-optimal) bound
\begin{equation*}
  t \ln \ln t \lesssim \mathbb{E}|X_t|^2 \lesssim t\ln t.
\end{equation*} 
In the recent work \cite{CHT22}, the following improved bound was obtained
\begin{equation*}
(\ln \ln t)^{-1-} \lesssim \frac{\mathbb{E}|X_t|^2}{t\sqrt{\ln t}} \lesssim (\ln \ln t)^{1+},
\end{equation*}
which essentially proves the conjectured $\sqrt{\ln t}$-super-diffusivity (up to a double-logarithmic correction).  

\smallskip

Compared to \cite{CHT22}, our approach follows a completely different homogenization-based strategy.  In addition to removing the double-logarithmic corrections, we are able to study the real time (as opposed to Laplace transform) asymptotics of the mean-square displacement directly using PDEs associated with the process.

\smallskip

The iterative approach introduced here does not rely heavily on the specific features of the GFF and, therefore, in principle, can be adapted to other problems.  For instance, in \cite{KO02}, Komorowski and Olla study the problem analogous to ours, except that the power spectrum of the field diverges algebraically at the origin.  We plan to apply our methods to that setting.


\section{Heuristics: an ODE for $\lambda_L$ as function of $L$} \label{s:heur}

In this section, we recall homogenization theory to 
motivate why the effective diffusivity $\lambda_L$ approximately satisfies
the following differential equation as a function of the IR cut-off scale $L\gg 1$:
\begin{align}\label{ao07}
\lambda_L\frac{\dif \lambda_L}{\dif L}\approx\frac{\epsilon^2}{2}\frac{\dif }{\dif L}\mathbb{E}\psi_L^2.
\end{align}
Since we may rewrite (\ref{ao07}) as 
$\frac{\dif }{\dif L}(\lambda_L^2-\epsilon^2\mathbb{E}\psi_L^2)\approx0$, 
it yields by integration from $L=1$
\begin{align*}
\lambda_L^2\approx1+\epsilon^2\mathbb{E}\psi_L^2,
\end{align*}
into which we insert (\ref{ao37}) to obtain Theorem \ref{t:main}. Note that \eqref{ao07} has a similar spirit to the heuristics from \cite{FFQSS85} which we present in Appendix \ref{sec:appendix_heuristic}.

\smallskip
Since the quadratic part of the stationary (and ergodic) matrix field $a_L$ is positive, 
\begin{align}\label{ao40}
\xi\cdot a_L\xi=|\xi|^2,
\end{align}
we are thus in the framework of the homogenization of (non-symmetric)
elliptic equations in divergence form.
The only subtlety is that while the stationary $a_L$ has Gaussian 
moments (however depending on $L$), it is not deterministically bounded from 
above\footnote{Identifying sufficient moment bounds on non-uniformly elliptic stationary 
and ergodic coefficient
fields is an active line of research, popularized in \cite{B11} and \cite{CD16}, 
see the recent \cite{BS22} for the sharpest results available.}. The first result we need, which is the starting point for the homogenization theory of divergence-form elliptic operators, concerns the existence of so-called correctors.

\begin{theorem} Given any vector $\xi \in \mathbb{R}^{2}$, there exists a unique stationary gradient field $\nabla u_{L}$ such that $\mathbb{E} \nabla u_{L} = \xi$ and for which the corresponding flux field $a_{L} \nabla u_{L}$ is divergence-free.  In particular, the random function $u_{L}$ is $a_{L}$-harmonic, meaning
	\begin{equation*}
		- \nabla \cdot a_{L} \nabla u_{L} = 0 \quad \text{in} \, \, \mathbb{R}^{2},
	\end{equation*}
and the gradient $\nabla u_{L}$ depends linearly on $\xi$. \end{theorem}

	\begin{proof} Observe that $\mathbb{E} |a_{L}(x)|^{p} < \infty$ for all $p < \infty$ and $\xi \cdot a_{L}(x) \xi = |\xi|^{2}$ for all $(x,\xi) \in \mathbb{R}^{d} \times \mathbb{R}^{d}$.  Thus, the existence and uniqueness of $\nabla u_{L}$ follows from \cite[Lemma 1]{BFO18}.  (Alternatively, the same question is treated in the specific setting of divergence-free drifts in \cite{Fe20}.) Linearity of $\nabla u_{L}$ in $\xi$ follows from uniqueness. \end{proof}

\smallskip

The effective diffusivity $\bar{a}_{L} \xi$ in direction $\xi$ is, by definition, the expectation of the flux $\bar{a}_{L} \xi = \mathbb{E}a_{L} \nabla u_{L}$ (so that $\bar{a}_{L}$ relates the averaged field to the averaged flux/current, like $a_{L}$ relates the microscopic field to the microscopic flux).  One typically writes $\nabla u_{L} = \xi + \nabla \phi_{L}$, such that the stationary mean-zero gradient $\nabla \phi_{L}$ satisfies
	\begin{align} \label{ao01}
		- \nabla \cdot a_{L} (\xi + \nabla \phi_{L}) = 0 \quad \text{with} \quad \mathbb{E} \nabla \phi_{L} = 0.
	\end{align}
Since the Ansatz amounts to $u_{L}(x) = \xi \cdot x + \phi_{L}(x)$, one speaks of $u_{L}$ as providing a harmonic coordinate function (which transforms the diffusion into a martingale, see Appendix \ref{sec:appendix_prelim} below), where $\phi_{L}$ ``corrects" the affine coordinate, thus the name corrector.  
Since $\nabla u_{L} = \xi + \nabla \phi_{L}$, the corrector determines the effective diffusivity $\bar{a}_{L}$ via the formula
	\begin{equation} \label{E: effective diffusivity with corrector}
		\bar{a}_{L} \xi = \mathbb{E}a_{L} (\xi + \nabla \phi_{L}).
	\end{equation}
The next lemma makes precise the fact that, by (statistical) isotropy of $\psi_{L}$, $\bar{a}_{L}$ is isotropic.

	\begin{lemma} \label{L: isotropy} There exists a number $\lambda_{L}$ such that $\bar{a}_{L} = \lambda_{L} \text{Id}$.  Further, for any unit vector $\xi \in \mathbb{R}^{2}$,
		\begin{align} \label{ao08}
	 		\lambda_{L} = \mathbb{E}|\xi + \nabla \phi_{L}|^{2} = 1 + \mathbb{E}|\nabla \phi_{L}|^{2} \geq 1.
		\end{align}
	\end{lemma}
	
The constant $\lambda_{L}$ is exactly the one that appeared earlier in \eqref{e:diffusivity-X_L}.  For more on this point, we refer the readers to Appendix \ref{sec:appendix_prelim}.
	
	\begin{proof} For the sake of precision, in this proof, we write $\nabla \phi_{L \xi}$, making the dependence on $\xi$ explicit.  For any rotation matrix $R$, the rotated field $\psi_{L}(R x)$ has the same law as $\psi_{L}(x)$.    It follows that, for any unit vector $\xi$, the joint law of the fields $(a_{L}(Rx), R^{*} \nabla \phi_{L \xi}(Rx))$ is the same as that of $(a_{L}(x),\nabla \phi_{L (R^{*}\xi)}(x))$.  From this, the identity \eqref{E: effective diffusivity with corrector}, and the fact that $J$ commutes with rotations, we deduce that 
		\begin{align*}
			\bar{a}_{L} (R^{*} \xi) = \mathbb{E}[a_{L}(Rx) R^{*} (\xi + \nabla \phi_{L \xi}(Rx))] = R^{*} (\bar{a}_{L} \xi).
		\end{align*}
	Thus there is a number $\lambda_{L}$ such that $\bar{a}_{L} = \lambda_{L} \text{Id}$.  
	
	Finally, given any unit vector $\xi$, since $-\nabla \cdot a_{L} (\xi + \nabla \phi_{L}) = 0$ and the quadratic form of $a_{L}$ is given by \eqref{ao40}, we compute
		\begin{align*}
			\lambda_{L} = \xi \cdot \bar{a}_{L} \xi = \mathbb{E} \xi \cdot a_{L}(\xi + \nabla \phi_{L \xi}) = \mathbb{E} (\xi + \nabla \phi_{L \xi}) \cdot a_{L} (\xi + \nabla \phi_{L \xi}) = \mathbb{E}|\xi + \nabla \phi_{L \xi}|^{2}.
		\end{align*}
	Further, since $|\xi| = 1$ and $\mathbb{E} \nabla \phi_{L \xi} = 0$, this is equivalent to $\lambda_{L} = 1 + \mathbb{E}|\nabla \phi_{L \xi}|^{2}$. \end{proof}

\smallskip
 
Motivated by \cite{SZ06}, we will now (heuristically) derive \eqref{ao07} 
on the level of finite differences, 
which amounts to studying the progressive homogenization with respect 
to increasing spatial scales. Toward that end, it will be useful to explicitly couple the fields $(\psi_{L})_{L \geq 1}$.  As in the introduction, let $b$ be divergence-free vector white noise with UV cut-off at scale one, that is, the stationary centered Gaussian vector field with correlation function $R$ determined by \eqref{E: spectrum of b}.  For any $L \geq 1$, let $\psi_{L}$ be the scalar field determined by $b$ through the formula
	\begin{align}
		\mathcal{F}\psi_{L}(k) = |k|^{-2} I(L^{-1} \leq |k| \leq 1) Jk \cdot \mathcal{F}b(k),
	\end{align}
where $\mathcal{F}$ denotes the (distributional) Fourier transform.  In view of the definition of $b$, $\psi_{L}$ is a stationary Gaussian field with correlation function $c_{L}$ determined by \eqref{ao39}.  The next lemma asserts that, as a process indexed by $L$, $\psi_{L}$ has independent increments and the corrector field $\nabla \phi_{L}$ is adapted to the associated filtration.  

	\begin{lemma} \label{L: independence lemma} The process $(\psi_{L})_{L \geq 1}$ is a martingale with independent increments.  In particular, for any $L \leq L_{+}$, the field $\psi_{L}$ is independent of the increment $\psi_{L}' := \psi_{L_{+}} - \psi_{L}$.  Furthermore, $\nabla \phi_{L}$ is measurable with respect to the $\sigma$-algebra generated by $\psi_{L}$. \end{lemma} 
	
		\begin{proof} Given any two scales $L \leq L_{+}$, the field $\psi_{L}$ and the increment $\psi_{L}'$ are uncorrelated since, for any two Schwarz functions $\zeta_{1}, \zeta_{2}$,
	\begin{align*}
		\mathbb{E} (\zeta_{1} * \psi_{L})(0) (\zeta_{2} * \psi_{L}')(0)  = \int_{\{L^{-1} \leq |k| \leq 1\} \cap \{L_{+}^{-1} \leq |k| \leq L^{-1}\}} \frac{ \dif k}{2 \pi} |k|^{-2} \mathcal{F}\zeta_{1}(k) \overline{\mathcal{F}\zeta_{2}(k)} = 0.
	\end{align*}
Thus, by Gaussianity, $\psi_{L}$ and $\psi_{L}'$ are independent. Finally, since the coefficient $a_{L}$ appearing in the PDE \eqref{ao01} is measurable with respect to the $\sigma$-algebra generated by $\psi_{L}$, it follows that the solution $\nabla \phi_{L}$ is also. \end{proof}
	
\smallskip
	
In the rest of this section, we explain how to exploit the martingale structure of $\psi_{L}$ to analyze the asymptotics of $\lambda_{L}$.  Comparing two scales $L\le L_+$, we regroup the difference
of the corresponding corrector equations, see (\ref{ao01}), as
\begin{align}\label{ao02}
-\nabla\cdot a_L\nabla(\phi_{L_+}-\phi_L)
=\epsilon\nabla\cdot(\psi_{L+}-\psi_L)J(\xi+\nabla\phi_{L+}).
\end{align}
We now carry out two approximations in (\ref{ao02}):
\begin{itemize}
\item {\sc Homogenization}. Since we expect $\phi_{L+}-\phi_L$
to mostly live on scales $\gtrsim L$, we replace the elliptic operator 
$-\nabla\cdot a_L\nabla$ on the {l.h.s.} of (\ref{ao02})
by its homogenization $-\nabla\cdot\lambda_L\nabla = -\lambda_{L} \Delta$. 
\item {\sc Linearization}. Since $\phi_L$ is of first order in $\epsilon$, 
the {r.h.s.} contribution $\epsilon\nabla\cdot(\psi_{L_+}-\psi_L)J\nabla\phi_{L_+}$ is of 
second order in $\epsilon$, which we neglect.
\end{itemize}
Hence we approximate the increment $\phi_{L_+}-\phi_L$ by the stationary Gaussian $\epsilon\phi_L'$ 
defined through
\begin{align}\label{ao03}
-\lambda_{L} \Delta\phi_L'=\nabla\cdot\psi_L'J\xi\quad\mbox{with}\quad
\mathbb{E}\phi_L'=0,
\end{align}
and where, as in Lemma \ref{L: independence lemma}, we have set for abbreviation
\begin{align}\label{ao26}
\psi_L':=\psi_{L_+}-\psi_L.
\end{align}
More precisely, we approximate
$\phi_{L_+}-\phi_L$ by what has come to be called the two-scale expansion of $\epsilon\phi_L'$ (see \cite{allaire} and the references therein),
which endows the large-scale $\epsilon\phi_L'$ with appropriate smaller-scale modulations, 
see also Section \ref{s:sket},
\begin{align}\label{ao17}
\phi_{L+}\approx\phi_L+\epsilon(1+\phi_{Li}\partial_i)\phi_L'.
\end{align}
In (\ref{ao17}), $\phi_{Li}$ denotes the solution
of (\ref{ao01}) with the generic $\xi$ 
replaced by the Cartesian unit vector $e_i$, $i=1,\cdots,d$;
where here and in the sequel we use Einstein's convention of implicit summation over
repeated indices. By linearity, we have
\begin{align}\label{as38}
\phi_L=\xi_i\phi_{Li}.
\end{align}
When passing to the gradient of
(\ref{ao17}), we ignore the second-order derivatives
of $\phi_L'$ (since the latter lives on scales $\ge L$), 
which amounts to
\begin{align}\label{ao09}
\nabla\phi_{L_+}\approx\nabla\phi_L+\epsilon\partial_i\phi_L'(e_i+\nabla\phi_{Li}).
\end{align}

\smallskip

We now come to the core probabilistic argument which leverages the martingale structure: Since $\psi_L'$ is independent
of $\psi_L$, and both are stationary, we not only have
\begin{align}\label{ao05}
\mathbb{E}\psi_{L_+}^2=\mathbb{E}\psi_L^2+\mathbb{E}{\psi_{L}'}^2,
\end{align}
but also the stationary $\nabla\phi_L'$ defined through (\ref{ao03}), being measurable with respect to the $\sigma$-algebra generated by $\psi_{L}'$, is 
independent of (the stationary) $\nabla\phi_L$. Hence when applying $\mathbb{E}|\cdot|^2$
to the {r.h.s.} of (\ref{ao09}), the mixed term $2\epsilon\mathbb{E}\partial_i\phi_L'\nabla\phi_L\cdot(e_i+\nabla\phi_{Li})$ factorizes to $2\epsilon\mathbb{E}\partial_i\phi_L'
\mathbb{E}\nabla\phi_L\cdot(e_i+\nabla\phi_{Li})$ and thus vanishes because the first
factor vanishes by the second item in (\ref{ao01}). 
Therefore the approximation (\ref{ao09}) leads to
\begin{align*}
\mathbb{E}|\xi+\nabla\phi_{L_+}|^2\approx\mathbb{E}|\xi+\nabla\phi_L|^2
+\epsilon^2\mathbb{E}(\partial_i\phi_L')\partial_j\phi_L'
\mathbb{E}(e_i+\nabla\phi_{Li})\cdot(e_j+\nabla\phi_{Lj}).
\end{align*} 
We now appeal to an isotropy argument: Since by (\ref{ao08}), $\mathbb{E}|\xi+\nabla\phi_L|^2$
is independent of the unit vector $\xi$, we obtain from (\ref{as38}) that\footnote{Here
$\delta_{ij}$ denotes the Kronecker symbol.}
$\mathbb{E}(e_i+\nabla\phi_{Li})\cdot(e_j+\nabla\phi_{Lj})$ 
$=\delta_{ij}\mathbb{E}|\xi+\nabla\phi_L|^2$, so that the above collapses to
\begin{align}\label{ao19}
\mathbb{E}|\xi+\nabla\phi_{L_+}|^2\approx\mathbb{E}|\xi+\nabla\phi_L|^2
+\epsilon^2\mathbb{E}|\nabla\phi_L'|^2
\mathbb{E}|\xi+\nabla\phi_{L}|^2.
\end{align}
According to (\ref{ao08}) this can be reformulated by
\begin{align}\label{ao06}
\lambda_{L_+}\approx\lambda_L(1+\epsilon^2\mathbb{E}|\nabla\phi_L'|^2).
\end{align}
In Lemma \ref{l:prime}, we shall give the (easy) argument for what amounts
to an equi-partition result
\begin{align}\label{ao34}
\lambda_L^2\mathbb{E}|\nabla\phi_L'|^2
=\frac{1}{2}\mathbb{E}{\psi_L'}^2.
\end{align}
Inserting (\ref{ao34}) and (\ref{ao05}) into (\ref{ao06}) yields
\begin{align}\label{ao12}
\lambda_{L_+}\approx\lambda_L\Big(1+\frac{\epsilon^2}{2\lambda_L^2}
\big(\mathbb{E}\psi_{L_+}^2-\mathbb{E}\psi_L^2\big)\Big),
\end{align}
which is the finite-difference version of (\ref{ao07}).

\smallskip

Two comments on this heuristics are in place: 
If we had ignored the two-scale correction in (\ref{ao09}), {i.e.}
if we had approximated $\nabla\phi_{L_+}$ just by $\nabla\phi_L+\epsilon\nabla\phi_L'$,
then we would have obtained $\lambda_L^2\frac{\dif \lambda_L}{\dif L}$ on the {l.h.s.} of (\ref{ao07}),
which would have been the wrong approximation.
If we had ignored $O(\epsilon^2)$-terms in the ``cumulative'' (\ref{ao01})
instead of the ``incremental'' (\ref{ao02}), we would have ended up with
$\lambda_L=1+\frac{\epsilon^2}{2}\mathbb{E}\psi_L^2$, which would have translated into  
$\frac{\dif \lambda_L}{\dif L}$ on the {l.h.s.} of (\ref{ao07}), again the wrong approximation
for $\epsilon^2\ln L\gtrsim 1$.


\section{Homogenization error estimate 
based on augmented corrector $(\phi_L,\sigma_L)$}\label{s:sket}

In this section, we introduce our general approach to estimating the homogenization
error, which also deals with the linearization error. A key object is 
the flux corrector $\sigma_L$, which complements the ``field'' corrector $\phi_L$,
see \cite[Lemma 1]{GNO20} for its introduction into the realm of (qualitative) 
stochastic homogenization.
In our case of two space dimensions, there is no gauge freedom for a vector potential, 
and both correctors can be conveniently characterized by
\begin{align}\label{ao18}
a_L(\xi+\nabla\phi_L)=\lambda_L\xi+J\nabla\sigma_L\quad\mbox{with}\quad
\mathbb{E}\nabla\phi_L=\mathbb{E}\nabla\sigma_L=0.
\end{align}
(Though the existence and uniqueness of $\nabla \sigma_{L}$ will not be needed in the proofs that follow, the result is recalled in Appendix \ref{sec:appendix_prelim} below.)  Note that by the second item in (\ref{ao95}), 
applying the divergence to (\ref{ao18}) yields (\ref{ao01}), and taking
the expectation of (\ref{ao18}) yields (\ref{ao08}) in its original form of
$\lambda_L\xi$ $=\mathbb{E}a_L(\xi+\nabla\phi_L)$. 
The merit of (\ref{ao18}) (with the generic $\xi$ replaced by the Cartesian $e_i$), 
is that it shows that the two-scale expansion connects $a_L$ to $\lambda_L$
on the level of the flux: For any function $\zeta$ we have by Leibniz' rule
\begin{align}\label{ao23}
a_L\nabla(1+\phi_{Li}\partial_i)\zeta
=\lambda_L\nabla\zeta+J\nabla\sigma_{Li}\partial_i\zeta
+(\phi_{Li}a_L-\sigma_{Li}J)\nabla\partial_i\zeta,
\end{align}
see for instance \cite[(11)]{JO22} for the straightforward calculation 
leading from (\ref{ao18}) to (\ref{ao23}).

\smallskip

While in our actual proof we will not quite proceed like this, we now sketch how
the augmented corrector serves in an estimate of the homogenization (and linearization) error.
Identity (\ref{ao23}) 
allows to derive a divergence-form equation for the approximation error in (\ref{ao17}):
Using Leibniz' rule in the same vein as in the passage from (\ref{ao18}) to (\ref{ao23})
to rewrite the linearization error 
$\psi_L'J\nabla\phi_{L_+}$ $=-\phi_{L_+}J\nabla\psi_L'+J\nabla\psi_L'\phi_{L_+}$,
using (\ref{ao23}) for $\zeta=\phi_L'$, applying the divergence,
appealing to (\ref{ao03}), and subtracting (\ref{ao02}) from the result, we obtain
\begin{align}\label{ao22}
-\nabla\cdot a_L\nabla\big(\phi_{L_+}-(\phi_L+\epsilon(1+\phi_{Li}\partial_i)\phi_L')\big)
=\epsilon\nabla\cdot\big((\phi_{Li}a_L-\sigma_{Li}J)\nabla\partial_i\phi_L'
-\phi_{L+}J\nabla\psi_L'\big).
\end{align}
The first {r.h.s.} term $(\phi_{Li}a_L-\sigma_{Li}J)\nabla\partial_i\phi_L'$
results from homogenization, the second term $\phi_{L+}J\nabla\psi_L'$ from linearization. 
The merit of representation (\ref{ao22}) is that in
$(\phi_{Li}a_L-\sigma_{Li}J)\nabla\partial_i\phi_L'$,
there is no gradient on the smaller-scale objects $(\phi_L,\sigma_L)$,
but two gradients on the large-scale object $\phi_L'$.

\smallskip

In view of (\ref{ao40}) and the energy estimate, 
formula (\ref{ao22}) allows to estimate the approximation error
in (\ref{ao09}) with respect to $\mathbb{E}^\frac{1}{2}|\cdot|^2$, 
which is the norm relevant for the argument (\ref{ao19}), see Lemma \ref{l:proxy}. 
It is estimated by 
$\epsilon\big((\phi_{Li}a_L-\sigma_{Li}J)\nabla\partial_i\phi_L'$
$-\phi_{L+}J\nabla\psi_L'\big)$ in the same norm. To leading order in $\epsilon$,
the latter expression reduces to 
$\epsilon\big((\phi_{Li}-\sigma_{Li}J)\nabla\partial_i\phi_L'$
$-\phi_{L}J\nabla\psi_L'\big)$. Momentarily ignoring $\sigma_L$, 
which in view of Lemma \ref{l:est} is qualitatively similar to $\phi_L$,
we obtain by the triangle inequality and the martingale structure
\begin{align*}
\epsilon\mathbb{E}^\frac{1}{2}\big|\phi_{Li}\nabla\partial_i\phi_L'-\phi_{L}J\nabla\psi_L'\big|^2
\le\epsilon\big({\big(\mathbb{E}\phi_{Li}\phi_{Lj}\mathbb{E}\nabla\partial_i\phi_L'\cdot
\nabla\partial_j\phi_L'}\big)^\frac{1}{2}
+\mathbb{E}^\frac{1}{2}\phi_{L}^2\mathbb{E}^\frac{1}{2}|\nabla\psi_L'|^2\big),
\end{align*}
which by a similar isotropy argument (see (\ref{orthogonality})) as before (\ref{ao19}) yields
\begin{align}\label{ao41}
\epsilon\mathbb{E}^\frac{1}{2}\big|\phi_{Li}\nabla\partial_i\phi_L'-\phi_{L}J\nabla\psi_L'\big|^2
\le\epsilon{\mathbb{E}^\frac{1}{2}\phi_{L}^2}
\big(\mathbb{E}^\frac{1}{2}
{|\nabla^2\phi_L'|^2}
+\mathbb{E}^\frac{1}{2}|\nabla\psi_L'|^2\big),
\end{align}
where here
and in the sequel, $|\nabla^2\phi_L'|$ denotes the Frobenius norm of the
Hessian matrix $\nabla^2\phi_L'$.
The second factor is (easily) estimated, see (\ref{ao54}) of Lemma \ref{l:prime};
because of $\lambda_L\ge 1$, see (\ref{ao08}), we thus may deduce from (\ref{ao41}) that
\begin{align*}
\epsilon^2\mathbb{E}\big|\phi_{Li}\nabla\partial_i\phi_L'-\phi_{L}J\nabla\psi_L'\big|^2
\lesssim\epsilon^2 L^{-2}\mathbb{E}\phi_{L}^2,
\end{align*}
where we use $A\lesssim B$ to denote $A\le C B$ for some constant $C$ that only depends on $\epsilon^2$, is uniform in $\epsilon^2$ when $\epsilon^2 \le 1$ but may worsen when $\epsilon^2 \ge 1$. Estimating $\mathbb{E}\phi_L^2$ seems out of reach -- in fact, we don't even know whether or not $\phi_L$ is stationary.
Hence in our proof, we will replace $\phi_L$ by a proxy $\tilde\phi_L$ that is stationary
by construction, see the next Section \ref{s:constr}, 
and then compare $\nabla\tilde\phi_L$ to the stationary $\nabla\phi_L$, see Lemma \ref{l:proxy}.


\section{Construction of a stationary proxy $(\tilde\phi_L,\tilde\sigma_L)$ for 
$(\phi_L,\sigma_L)$}\label{s:constr}

We make the Ansatz of a geometrically growing partitioning, i.e., we start at $L = 1$ and proceed recursively, passing from the current scale $L$ to the next scale $L_{+}$ via the rule
\begin{align}\label{ao76}
L_+=ML
\end{align}
for a sufficiently large factor $M$.  This restricts
$L$ to be an integer power of $M$, but this will be no restriction for us, see the proof of Theorem \ref{t:main}.  As in the discussion above, let $\psi_{L}' = \psi_{L_{+}} - \psi_{L}$.  

\smallskip

The construction of the proxy corrector $(\tilde{\phi}_{L},\tilde{\sigma}_{L})$ is recursive.  In addition to the correctors, the recursion tracks a proxy diffusivity $\tilde{\lambda}_{L}$ as well.  Beginning at scale $L = 1$, we set 
\begin{align*}
\tilde{\phi}_{L = 1} = 0, \quad \tilde{\sigma}_{L = 1} = 0, \quad \tilde{\lambda}_{L = 1} = 1.
\end{align*}

\smallskip

Given the proxy $(\tilde{\phi}_{L},\tilde{\sigma}_{L})$ and diffusivity $\tilde{\lambda}_{L}$ at some scale $L$, to proceed to scale $L_{+} = ML$, we begin by invoking the Helmholtz decomposition of the incremental field $\psi_{L}'J$: For any $\xi \in \mathbb{R}^{d}$, let $\phi_{L}'$ and $\sigma_{L}'$ be the stationary Gaussian fields characterized by 
\begin{align}\label{ao27}
\tilde\lambda_L\nabla\phi_L'+\psi_L'J\xi=J\nabla\sigma_L'\quad\mbox{and}
\quad\mathbb{E}\phi_L'=\mathbb{E}\sigma_L'=0.
\end{align}
Applying successively the divergence operator $\nabla \cdot$ and the curl operator $\nabla \cdot J$ to the expression above and using the second item in (\ref{ao95}), we observe that $\phi_{L}'$ and $\sigma_{L}'$ are equivalently characterized as the mean-zero stationary solutions of the equations
\begin{align}
-\tilde{\lambda}_{L} \Delta \phi_{L}' &= \nabla \cdot \psi_{L}' J \xi, \label{as03} \\
- \Delta \sigma_{L}' &= - \nabla \cdot \psi_{L}' \xi. \label{as03alt}
\end{align}
In particular, $\phi_{L}'$ is defined consistently with (\ref{ao03}) in Section \ref{s:heur}.  Note that both $\phi_{L}'$ and $\sigma_{L}'$ are stationary Gaussian fields, measurable with respect to the $\sigma$-algebra generated by the increment $\psi_{L}'$, and hence independent of $\psi_{L}$. The following lemma collects elementary but crucial properties of 
the incremental $(\phi_L',\sigma_L')$, the first one being (stochastic) isotropy,
namely a certain covariance in law under the change of variables $x\mapsto Rx$
for some orthogonal $R$.  The definition of the proxy correctors $(\tilde{\phi}_{L},\tilde{\sigma}_{L})$ resumes after the proof of the lemma.

\begin{lemma}\label{l:prime}
For a fixed but arbitrary orthogonal $R$, 
let $(\phi_L'^*,\sigma_L'^*)$ denote the solution of (\ref{ao27})
with $\xi$ replaced by $R^*\xi$. Then 
\begin{align}\label{as12}
(\phi_L'^*,\sigma_L'^*,\psi_L')=_{\mbox{law}}(\phi_L',({\rm det R})\sigma_L',
({\rm det R})\psi_L')(R\cdot).
\end{align}
Moreover, we have
\begin{align}
\tilde\lambda_L^2\mathbb{E}|\nabla\phi_L'|^2
&=\tilde\lambda_L\mathbb{E}\psi_L'J\nabla\phi_L'\cdot\xi
=\mathbb{E}|\nabla\sigma_L'|^2
=\frac{1}{2}\mathbb{E}{\psi_L'}^2
=\frac{1}{2}\ln M,\label{ao34bis}\\
\tilde\lambda_L^2\mathbb{E}|\nabla^2\phi_L'|^2
&\le\mathbb{E}|\nabla\psi_L'|^2
\le\frac{1}{2L^2},\label{ao54}\\
\tilde\lambda_L^2\mathbb{E}{\phi_L'}^2
&=\mathbb{E}{\sigma_L'}^2
\le \frac{L_+^2}{4}.\label{ao45}
\end{align}
\end{lemma}

{\sc Proof of Lemma \ref{l:prime}}. 
We start with (\ref{as12}) and note that because of 
\begin{align}\label{as15}
\nabla\zeta^*=R^*\nabla\zeta(R\cdot)\quad\mbox{provided}\quad\zeta^*=\zeta(R\cdot)
\quad\mbox{ and }\quad 
JR^*=({\rm det R})R^*J,
\end{align}
the triplet $(\phi_L',({\rm det R})\sigma_L',({\rm det R})\psi_L')(R\cdot)$
satisfies (\ref{ao27}) with $\xi$ replaced by $R^*\xi$. 
Since by isotropy (and Gaussianity), 
$({\rm det R})\psi_L'(R\cdot)$ has the same law as $\psi_L'$, the statement follows.

\smallskip

We note that the stationary centered fields
$\nabla\phi_L'$ and $J\nabla\sigma_L'$ are orthogonal with respect to $L^2$ in probability,
as can be seen by appealing to the second item in (\ref{ao95}) to obtain
$(\nabla\phi_L')\cdot J\nabla\sigma_L'$ =$\nabla\cdot\phi_L'J\nabla\sigma_L'$,
and then applying $\mathbb{E}$ so that by stationarity
$\mathbb{E}(\nabla\phi_L')\cdot J\nabla\sigma_L'=0$. Hence it follows from 
(\ref{ao27}) that $(-\tilde\lambda_L\nabla\phi_L',J\nabla\sigma_L')$ provides
an orthogonal decomposition of $\psi_L'J\xi$, which implies
\begin{align*}
\tilde\lambda_L^2\mathbb{E}|\nabla\phi_L'|^2+\mathbb{E}|\nabla\sigma_L'|^2
=\mathbb{E}{\psi_L'}^2,\quad
\tilde\lambda_L\mathbb{E}\psi_L'J\nabla\phi_L'\cdot\xi
=-\tilde\lambda_L\mathbb{E}\nabla\phi_L'\cdot\psi_L'J\xi
=\tilde\lambda_L^2\mathbb{E}|\nabla\phi_L'|^2.
\end{align*}
Comparing \eqref{as03} and \eqref{as03alt} (and using the same notation as in the proof of Lemma \ref{L: isotropy} to make the dependence on $\xi$ explicit), we notice that $\sigma_{L \xi}'  = \tilde\lambda_L\phi_{L (J\xi)}'$ for any $\xi$. Since by (\ref{as12}),
$\tilde\lambda_L^2\mathbb{E}|\nabla\phi_L'|^2$ is independent of
the unit vector $\xi$, we thus have
\begin{align}\label{as02}
\tilde\lambda_L^2\mathbb{E}|\nabla\phi_L'|^2=\mathbb{E}|\nabla\sigma_L'|^2.
\end{align}
The combination of these identities yields all but
the last identity in (\ref{ao34bis}), which follows from 
(\ref{ao37}), (\ref{ao26}), (\ref{ao05}),
and (\ref{ao76}).

\smallskip

For (\ref{ao54}), we note that based on the stationarity of $\phi_L'$ (and thus
also its derivatives) we have (using Einstein's summation convention as before)
\begin{align*}
\mathbb{E}(\partial_i\partial_j\phi_L')(\partial_i\partial_j\phi_L')
=\mathbb{E}(\partial_i\partial_i\phi_L')(\partial_j\partial_j\phi_L')
\end{align*}
which amounts to $\mathbb{E}|\nabla^2\phi_L'|^2$ $=\mathbb{E}(\Delta\phi_L')^2$. Therefore, the first inequality in (\ref{ao54}) follows from (\ref{as03}).
Now if $c_L'$ denotes the covariance function of the stationary $\psi_L'$,
we have 
\begin{align*}
\mathbb{E}|\nabla\psi_L'|^2=-\Delta c_L'(0)=\int\frac{\dif k}{2\pi}|k|^2{\mathcal F}c_L'.
\end{align*}
Since in view of (\ref{ao26}) we have $c_{L_+}=c_L+c_L'$, 
we obtain from (\ref{ao39}) 
\begin{align*}
\mathbb{E}|\nabla\psi_L'|^2
=\int_{\{L_+^{-1}\le|k|\le L^{-1}\}}\frac{\dif k}{2\pi}\le\frac{1}{2L^2}.
\end{align*}

\smallskip

We finally turn to (\ref{ao45}). The first identity follows for the same reason
as (\ref{as02}). Again, if $\bar c_L'$ is the covariance function of the
stationary $\phi_L'$, we have
\begin{align}\label{as08}
\mathbb{E}{\phi_L'}^2=\bar c_L'(0)=\int\frac{\dif k}{2\pi}{\mathcal F}\bar c_L'.
\end{align}
In view of (\ref{as03}), the covariance function $\bar c_L'$ is related 
to the covariance function $c_L'$ of $\psi_L'$ via
\begin{align*}
\tilde\lambda_L^2\Delta^2\bar c_L'=-J\xi\cdot \nabla^2c_L' J\xi,
\end{align*}
which on the level of the Fourier transform amounts to
\begin{align*}
\tilde\lambda_L^2|k|^4{\mathcal F}\bar c_L'(k)=(k\cdot J\xi)^2{\mathcal F}c_L'(k)
\stackrel{(\ref{ao39})}{=}I(L_{+}^{-1}\le|k|\le L^{-1})(k\cdot J\xi)^2|k|^{-2}.
\end{align*}
Inserting this into (\ref{as08}) we obtain 
\begin{align*}
\tilde\lambda_L^2\mathbb{E}{\phi_L'}^2
=\int_{\{L_{+}^{-1}\le|k|\le L^{-1}\}}\frac{\dif k}{2\pi}(k\cdot J\xi)^2|k|^{-6}
\le\frac{L_+^2}{2},
\end{align*}
which yields the remaining part of (\ref{ao45}).
\qed

\smallskip

Next, we define $\tilde\phi_{ L_{+} }$ in line with (\ref{ao17}):
\begin{align}
\tilde\phi_{L_+}&:=\tilde\phi_L+\epsilon(1+\tilde\phi_{Li}\partial_i)\phi_L'
\quad\mbox{and}\quad\tilde\phi_{L=1}\equiv0,\label{ao31}
\end{align}
where, as above, $\tilde\phi_{Li}$ is the version of $\tilde\phi_L$ with 
the generic unit vector $\xi$ replaced by the Cartesian $e_i$.
The recursive definition of the stationary proxy $\tilde\sigma_L$ of the flux corrector $\sigma_L$
has a similar structure, but appears less intuitive:
\begin{align}
\tilde\sigma_{L_+}&:=\tilde\sigma_L
+\epsilon\big(\sigma_L'+\tilde\phi_{Li}\psi_L'(\xi_i+\epsilon\partial_i\phi_L')
+\tilde\sigma_{Li}\partial_i\phi_L'\big)\quad\mbox{and}\quad
\tilde\sigma_{L=1}\equiv0;\label{ao32}
\end{align}
its merit will become apparent in the proof of Lemma \ref{l:induct} around (\ref{ao71}). Since $(\phi_L',\sigma_L')$ is obviously linear in $\xi$, 
linearity propagates through (\ref{ao31}) and (\ref{ao32}),
so that also $(\tilde\phi_L,\tilde\sigma_L)$ are linear in $\xi$, which gives the consistency
\begin{align}\label{as10}
\tilde\phi_L=\xi_i\tilde\phi_{Li}\quad\mbox{and}\quad
\tilde\sigma_L=\xi_i\tilde\sigma_{Li}.
\end{align}
%

\smallskip

In passing from the true corrector $(\phi_{L},\sigma_{L})$ to the proxy, the price we pay is that for the proxy $(\tilde\phi_L,\tilde\sigma_L)$, 
the flux decomposition (\ref{ao18}) holds only approximately:
\begin{align}\label{ao18bis}
a_L(\xi+\nabla\tilde\phi_L)=\tilde\lambda_L\xi+J\nabla\tilde\sigma_L+f_L
\quad\mbox{with}\quad \mathbb{E}f_L=0,
\end{align}
where the stationary centered vector-field $f_L$ captures the error. That constant deterministic vector is of the form $\tilde\lambda_L\xi$
which is a consequence of isotropy, meaning that the stationary vector field
$a_L(\xi+\nabla\tilde\phi_L)-J\nabla\tilde\sigma_L$ transforms according to (\ref{as16}): Indeed, on the one hand, 
this vector field,
and thus its expectation, is linear in $\xi$ as a consequence of (\ref{as10}). 
On the other hand, when we replace $\xi$ by $R^*\xi$, 
this expectation transforms in the same way. Hence the expectation has to be of the form 
$\tilde\lambda_L\xi$ for some (deterministic) scalar $\tilde\lambda_L$.
As a consequence, 
$(f_L,\tilde\lambda_L)$ is implicitly determined through (\ref{ao18bis}),
which is important for our inductive definition, since $\tilde\lambda_L$ 
feeds into the definition of the incremental objects $(\phi_L',\sigma_L')$, see (\ref{ao27}).
However, also $(f_L,\tilde\lambda_L)$ can be recovered inductively, 
which will be crucial for its estimate:

\begin{lemma}\label{l:induct}
\begin{align}
\tilde\lambda_{L_+}&=\tilde\lambda_L
\big(1+\frac{\epsilon^2\ln M}{2\tilde\lambda_L^2}\big)
\quad\mbox{and}\quad\tilde\lambda_{L=1}\equiv1,\label{ao30}\\
f_{L_+}&=f_L+\epsilon^2(1-\mathbb{E})\psi_L'J\nabla\phi_L'
+\epsilon\big(-\tilde\phi_{Li}(\xi_i+\epsilon\partial_i\phi_L')J\nabla\psi_L'\nonumber\\
&+(\tilde\phi_{Li}a_{L}-\tilde\sigma_{Li}J)\nabla\partial_i\phi_L'
+\partial_i\phi_L' f_{Li}\big)\quad\mbox{and}\quad f_{L=1}\equiv0,\label{ao55}
\end{align}
where, in line with our earlier convention, 
$(\tilde\sigma_{Li},f_{Li})$ denotes $(\tilde\sigma_L,f_L)$ 
with $\xi$ replaced by $e_i$.
\end{lemma}

Note that (\ref{ao30}) is consistent with (\ref{ao12}). We also note that the terms on the {r.h.s.} of (\ref{ao22}) 
essentially reappear on the {r.h.s.} of (\ref{ao55}). 

\smallskip

{\sc Proof of Lemma \ref{l:induct}}.
We note that under the perturbation (\ref{ao18bis}) of (\ref{ao18}),
the crucial identity (\ref{ao23}) turns into
\begin{align}\label{ao24}
a_L\nabla(1+\tilde\phi_{Li}\partial_i)\zeta
=\tilde\lambda_L\nabla\zeta+(\tilde\phi_{Li}a_L-\tilde\sigma_{Li}J)\nabla\partial_i\zeta
+J\nabla\tilde\sigma_{Li}\partial_i\zeta+\partial_i\zeta f_{Li}.
\end{align}
By definitions (\ref{ao25}) and (\ref{ao26}) we have 
\begin{align}\label{as30}
a_{L+}-a_L=\epsilon\psi_L'J.
\end{align}
Using the identity (\ref{ao24}) with $\zeta=\phi_L'$, 
we obtain for the flux increment
\begin{align*}
\lefteqn{a_{L_+}(\xi+\nabla\tilde\phi_{L_+})- a_{L}(\xi+\nabla\tilde\phi_{L})
\stackrel{(\ref{ao31}),(\ref{as30})}{=}\epsilon\big(\psi_L'J(\xi+\nabla\tilde\phi_{L_+})
+a_L\nabla(1+\tilde\phi_{Li}\partial_i)\phi_L'\big)}\nonumber\\
&=\epsilon\big(\psi_L'J(\xi+\nabla\tilde\phi_{L_+})
+\tilde\lambda_L\nabla\phi_L'
+(\tilde\phi_{Li}a_L-\tilde\sigma_{Li}J)\nabla\partial_i\phi_L'
+J\nabla\tilde\sigma_{Li}\partial_i\phi_L'+\partial_i\phi_L' f_{Li}\big).
\end{align*}
Next, we substitute $\psi_L'J\xi+\tilde\lambda_L\nabla\phi_L'$ by 
$J\nabla\sigma_L'$ according to (\ref{ao27}).
In order to eliminate $\tilde\phi_{L_+}$ on the {r.h.s.} we use the linearity of $\tilde\phi_L$
in $\xi$ to rewrite (\ref{ao31}) as 
\begin{align}\label{as31}
\tilde\phi_{L_+}=(\xi_i+\epsilon\partial_i\phi_L')\tilde\phi_{Li}+\epsilon\phi_L',
\end{align}
to which we apply $\psi_L'J\nabla$:
\begin{align*}
\lefteqn{a_{L_+}(\xi+\nabla\tilde\phi_{L_+})-a_{L}(\xi+\nabla\tilde\phi_{L})
=\epsilon^2\psi_L'J\nabla\phi_L'}\nonumber\\
&+\epsilon\Big(J\nabla\sigma_L'
+\psi_L'{J\nabla}(\xi_i+\epsilon\partial_i\phi_L')\tilde\phi_{Li}
+(\tilde\phi_{Li}a_{L}-\tilde\sigma_{Li}J)\nabla\partial_i\phi_L'
+J\nabla\tilde\sigma_{Li}\partial_i\phi_L'+\partial_i\phi_L' f_{Li}\Big).
\end{align*}
Clearly, it is advantageous to put gradients on 
large-scale incremental quantities like $\psi_L'$ 
rather than on the smaller-scale cumulative quantities like $\tilde\phi_L$.
Hence we write 
\begin{align*}
\psi_L'{J\nabla}(\xi_i+\epsilon\partial_i\phi_L')\tilde\phi_{Li}
=J\nabla\tilde\phi_{Li}\psi_L'(\xi_i+\epsilon\partial_i\phi_L')
-\tilde\phi_{Li}(\xi_i+\epsilon\partial_i\phi_L')J\nabla\psi_L'
\end{align*}
to obtain
\begin{align}\label{ao71}
a_{L_+}(\xi+\nabla\tilde\phi_{L_+})&-
a_{L}(\xi+\nabla\tilde\phi_{L})
=\epsilon J\nabla\big(\tilde\phi_{Li}\psi_L'(\xi_i+\epsilon\partial_i\phi_L')
+\sigma_L'+\tilde\sigma_{Li}\partial_i\phi_L'\big)+\epsilon^2\psi_L'J\nabla\phi_L'\nonumber\\
&+\epsilon\Big(-\tilde\phi_{Li}(\xi_i+\epsilon\partial_i\phi_L')J\nabla\psi_L'
+(\tilde\phi_{Li}a_{L}-\tilde\sigma_{Li}J)\nabla\partial_i\phi_L'
+\partial_i\phi_L' f_{Li}\Big).
\end{align}
This motivates the inductive
definition (\ref{ao32}) of $\tilde\sigma_L$, 
which allows us to rewrite (\ref{ao71}) more compactly as
\begin{align*}
\lefteqn{\big(a_{L_+}(\xi+\nabla\tilde\phi_{L_+})-J\nabla\sigma_{L_+}\big)-
\big(a_{L}(\xi+\nabla\tilde\phi_{L})-J\nabla\sigma_{L}\big)=\epsilon^2\psi_L'J\nabla\phi_L'}
\nonumber\\
&+\epsilon\Big(-\tilde\phi_{Li}(\xi_i+\epsilon\partial_i\phi_L')J\nabla\psi_L'
+(\tilde\phi_{Li}a_{L}-\tilde\sigma_{Li}J)\nabla\partial_i\phi_L'
+\partial_i\phi_L' f_{Li}\Big).
\end{align*}
By definition (\ref{ao18bis}) of $(f_L,\tilde\lambda_L)$ this assumes the form
\begin{align}\label{ao72}
\lefteqn{(\tilde\lambda_{L_+}-\tilde\lambda_L)\xi+f_{L_+}-f_L=\epsilon^2\psi_L'J\nabla\phi_L'}
\nonumber\\
&+\epsilon\Big(-\tilde\phi_{Li}(\xi_i+\epsilon\partial_i\phi_L')J\nabla\psi_L'
+(\tilde\phi_{Li}a_{L}-\tilde\sigma_{Li}J)\nabla\partial_i\phi_L'
+\partial_i\phi_L' f_{Li}\Big).
\end{align}

\smallskip

Applying $\mathbb{E}$ to (\ref{ao72}), appealing to the martingale structure,
to the second item in (\ref{ao18bis}), and to $\mathbb{E}\nabla\psi_L'=\mathbb{E}\nabla\partial_i\phi_L'=0$, we obtain 
\begin{align}\label{ao73}
(\tilde\lambda_{L_+}-\tilde\lambda_L)\xi=\epsilon^2\mathbb{E}\psi_L'J\nabla\phi_L'
{-\epsilon^2\mathbb{E}\tilde\phi_{Li}\mathbb{E}\partial_i\phi_L'J\nabla\psi_L'}
\stackrel{(\ref{ao54ter})}{=}
\epsilon^2\mathbb{E}\psi_L'J\nabla\phi_L'.
\end{align}
Applying $\xi\cdot$ to (\ref{ao73}) and by (\ref{ao34bis})
we get (\ref{ao30}).
Subtracting (\ref{ao73}) from (\ref{ao72}), we also obtain (\ref{ao55}).\qed

\smallskip

The following lemma shows that for $(\tilde\lambda_L,\nabla\tilde\phi_L)$ to be a good proxy for $(\lambda_L,\nabla\phi_L)$, we need to control $\mathbb{E}|f_L|^2$.

\begin{lemma}\label{l:proxy}
\begin{align}
&0\le\tilde\lambda_{L}^2-(1+\epsilon^2\ln L)\lesssim (\epsilon^2\ln M)
\big(\epsilon^2\ln M
+\ln(1+\epsilon^2\ln L)\big) \lesssim\epsilon^4(\ln M)\ln L,\label{ao79}\\
&\mathbb{E}|\nabla\phi_L-\nabla\tilde\phi_L|^2\le\mathbb{E}|f_L|^2,\label{ao70}\\
&|\lambda_L-\tilde\lambda_L|\lesssim\big(\tilde\lambda_L\mathbb{E}|f_L|^2\big)^\frac{1}{2}
+\mathbb{E}|f_L|^2.
\label{ao80}
\end{align}
\end{lemma}

We recall that the notation $A \lesssim B$ means that there is a universal constant $C$ (in particular, independent of $\epsilon$, $L$, and $M$) such that $A \leq C B$.

\smallskip

{\sc Proof of Lemma \ref{l:proxy}}. 
We start with (\ref{ao79}), noting that by (\ref{ao76}), 
identity (\ref{ao30}) can be reformulated as
\begin{align}\label{ao84}
\tilde\lambda_{L_+}^2-\epsilon^2\ln L_+
=\tilde\lambda_L^2-\epsilon^2\ln L
+\big(\frac{\epsilon^2\ln M}{2\tilde\lambda_L}\big)^2.
\end{align}
Because of $\tilde\lambda_{L=1}=1$  we obtain by iteration 
\begin{align*}
\tilde\lambda_L^2\ge1+\epsilon^2\ln L,
\end{align*}
which, in combination with (\ref{ao84}), implies
\begin{align*}
0\le(\tilde\lambda_{L_+}^2-\epsilon^2\ln L_+)-(\tilde\lambda_L^2-\epsilon^2\ln L)
\le\frac{\epsilon^4\ln^2M}{1+ \epsilon^2\ln L}.
\end{align*}
By our geometric partitioning (\ref{ao76}) we obtain by iteration in $\frac{\ln L}{\ln M}$ steps 
\begin{align*}
0\le\tilde\lambda_{L}^2-\epsilon^2\ln L
\le\epsilon^2\ln M\sum_{n=0}^{\frac{\ln L}{\ln M}-1}
\frac{\epsilon^2\ln M}{1+ n\epsilon^2\ln M}.
\end{align*}
We compare the sum to an integral,
\begin{align}\label{as36}
\sum_{n=0}^{\frac{\ln L}{\ln M}-1}
\frac{\epsilon^2\ln M}{1+ n\epsilon^2\ln M}
\le\epsilon^2\ln M+\int_0^{\frac{\ln L}{\ln M}}\dif n\frac{\epsilon^2\ln M}{1+ n\epsilon^2\ln M},
\end{align}
and compute the integral,
\begin{align}\label{as37}
\int_0^{\frac{\ln L}{\ln M}}\dif n\frac{\epsilon^2\ln M}{1+ n\epsilon^2\ln M}
=\ln(1+\epsilon^2\ln L),
\end{align}
which yields (\ref{ao79}).

\smallskip

Applying $\nabla\cdot$ to (\ref{ao18bis}), appealing to the second item in (\ref{ao95}), 
and taking the difference with (\ref{ao01}) we obtain
\begin{align*}
-\nabla\cdot a_L(\nabla\phi_L-\nabla\tilde\phi_L)=\nabla\cdot f_L.
\end{align*}
Since $\nabla\phi_L-\nabla\tilde\phi_L$ is stationary and of vanishing expectation, this implies
\begin{align*}
\mathbb{E}(\nabla\phi_L-\nabla\tilde\phi_L)\cdot a_L(\nabla\phi_L-\nabla\tilde\phi_L)=
-\mathbb{E}(\nabla\phi_L-\nabla\tilde\phi_L)\cdot f_L
\end{align*}
According to (\ref{ao40}), the {l.h.s.} simplifies to
$\mathbb{E}|\nabla\phi_L-\nabla\tilde\phi_L|^2$, so that (\ref{ao70}) follows
from the Cauchy-Schwarz inequality.

\smallskip

We finally turn to (\ref{ao80}) and to this purpose multiply (\ref{ao18bis}) by
$\xi+\nabla\tilde\phi_L$; by (\ref{ao40}) and the second item in (\ref{ao95}), this turns into
\begin{align*}
|\xi+\nabla\tilde\phi_L|^2=\tilde\lambda_L(1+\xi\cdot\nabla\tilde\phi_L)
-\nabla\cdot\tilde\sigma_L J(\xi+\nabla\tilde\phi_L)
+f_L\cdot(\xi+\nabla\tilde\phi_L).
\end{align*}
Applying $\mathbb{E}$, which commutes with $\nabla\cdot$, we obtain by stationarity
and by centeredness of $f_L$
\begin{align}\label{ao96}
1+\mathbb{E}|\nabla\tilde\phi_L|^2
=\mathbb{E}|\xi+\nabla\tilde\phi_L|^2=\tilde\lambda_L+\mathbb{E}f_L
\cdot\nabla\tilde\phi_L.
\end{align}
By Young's inequality, (\ref{ao96}) implies in particular
\begin{align}\label{ao97}
\mathbb{E}|\nabla\tilde\phi_L|^2\lesssim\tilde\lambda_L+\mathbb{E}|f_L|^2.
\end{align}

For \eqref{ao80}, we subtract (\ref{ao08}) from (\ref{ao96}):
\begin{align*}
2\mathbb{E}(\nabla\tilde\phi_L-\nabla\phi_L)\cdot\nabla\tilde\phi_L
-\mathbb{E}|\nabla\tilde\phi_L-\nabla\phi_L|^2
=\tilde\lambda_L-\lambda_L+\mathbb{E}f_L\cdot\nabla\tilde\phi_L.
\end{align*}
By Young's inequality, this implies (\ref{ao80}):
\begin{align*}
|\tilde\lambda_L-\lambda_L|
&\lesssim\mathbb{E}|\nabla\tilde\phi_L-\nabla\phi_L|^2
+\big((\mathbb{E}|f_L|^2+\mathbb{E}|\nabla\tilde\phi_L-\nabla\phi_L|^2)
\mathbb{E}|\nabla\tilde\phi_L|^2\big)^\frac{1}{2}\nonumber\\
&\stackrel{(\ref{ao70})}{\lesssim}\mathbb{E}|f_L|^2
+\big(\mathbb{E}|f_L|^2\mathbb{E}|\nabla\tilde\phi_L|^2\big)^\frac{1}{2}
\stackrel{(\ref{ao97})}{\lesssim}\mathbb{E}|f_L|^2
+\big(\tilde\lambda_L\mathbb{E}|f_L|^2\big)^\frac{1}{2}.
\end{align*}
\qed



\section{Estimate of the stationary proxy $(\tilde\phi_L,\tilde\sigma_L)$ } \label{s:approxcorr}

The following lemma establishes that $(\tilde\phi_L,\tilde\sigma_L)$ is
isotropic, centered, and of $O(\epsilon L)$, at least on the level of
second and fourth moments. We note that by (\ref{ao25}) 
the {r.h.s.}  of (\ref{ao55}) contains the term
$\epsilon\tilde\phi_{Li}\psi_L\nabla\partial_i\phi_L'$, 
which is quadratic (instead of just linear)
in the smaller-scale objects $(\tilde\phi_L,\tilde\sigma_L,\psi_L)$.
This necessitates the fourth-moment control (\ref{as17bis}).

\begin{lemma}\label{l:est} 
We have isotropy in the same sense as in Lemma \ref{l:prime}: 
For fixed but arbitrary orthogonal $R$,
let $(\tilde\phi_L^*,\tilde\sigma_L^*)$ denote the outcome of the iteration
(\ref{ao31}) and (\ref{ao32}) with $\xi$ replaced by $R^*\xi$. 
Then 
\begin{align}\label{as11}
(\tilde\phi_L^*,\tilde\sigma_L^*,\psi_L)=_{\mbox{law}}
\big(\tilde\phi_L,({\rm det}R)\tilde\sigma_L,
({\rm det} R)\psi_L\big){(R\cdot)}.
\end{align}
As a consequence,
\begin{align}\label{orthogonality}
\mathbb{E}\tilde{\phi}_{Li}\tilde{\phi}_{Lj}=\delta_{ij}\mathbb{E}\tilde{\phi}_{L}^{2}
\quad\mbox{and}\quad
\mathbb{E}\tilde{\sigma}_{Li}\tilde{\sigma}_{Lj}=\delta_{ij}\mathbb{E}\tilde{\sigma}_{L}^{2}.
\end{align}
The pair $(\tilde\phi_L,\tilde\sigma_L)$ is centered, that is,
\begin{align}
\mathbb{E}\tilde\phi_L=\mathbb{E}\tilde\sigma_L&=0\label{ao54ter}
\end{align}
and for $M\gg1$, it satisfies the estimates
\begin{align}
\tilde\lambda_L^2\mathbb{E}\tilde\phi_L^2+\mathbb{E}\tilde\sigma_{L}^2
&\lesssim \epsilon^2 L^2,\label{ao64}\\
\tilde\lambda_L^4\mathbb{E}\tilde\phi_L^4
&\lesssim \epsilon^4 L^4.\label{as17bis}
\end{align}
\end{lemma}

\smallskip

In the statement above and also in Lemma \ref{s:error} below, by $M \gg 1$ we mean that there is an $M_{0} < \infty$ (independent of $L$ and $\epsilon$) such that the statement holds provided $M \geq M_{0}$.

\smallskip

{\sc Proof of Lemma \ref{l:est}}. We establish (\ref{as11}) inductively:
It follows from the induction hypothesis in conjunction with the 
linearity (\ref{as10}) that directional derivatives transform as
\begin{align*}
\tilde\phi_{Li}^*\partial_i\zeta^*=\tilde\phi_{Lj}\partial_j\zeta(R\cdot)\quad\mbox{and}\quad
\tilde\sigma_{Li}^*\partial_i\zeta^*=({\rm det}R)\tilde\sigma_{Lj}\partial_j\zeta(R\cdot)
\quad\mbox{provided}\quad\zeta^*=\zeta(R\cdot).
\end{align*}
Together with the incremental transformation behavior (\ref{as12}) we learn
that (\ref{ao31}) and (\ref{ao32}) preserve the transformation behavior (\ref{as11}). 
In case of (\ref{ao32}) we made use of the fact that the {r.h.s.} is linear in the quantities 
$(\psi_L',\sigma_L',\tilde\sigma_L)$ for which the transformation behavior
involves the factor ${\rm det}R$. 

\smallskip

Next, we turn to (\ref{orthogonality}).  
On the one hand by (\ref{as11}), $\mathbb{E}\tilde\phi_{L}^2$ 
does not depend on the unit vector $\xi$; 
on the other hand, by (\ref{as10}), it can be expressed as
$\xi_i\xi_j\mathbb{E}\tilde\phi_{Li}\tilde\phi_{Lj}$; hence we must have
\begin{align*}
\mathbb{E}\tilde\phi_{Li}\tilde\phi_{Lj}=\delta_{ij}\mathbb{E}\tilde\phi_L^2
\quad\mbox{and likewise}\quad
\mathbb{E}\tilde\sigma_{Li}\tilde\sigma_{Lj}=\delta_{ij}\mathbb{E}\tilde\sigma_L^2.
\end{align*}

\smallskip

We now turn to (\ref{ao54ter}).
 Applying $\mathbb{E}$ to (\ref{ao31}), we learn from the martingale structure
\begin{align*}
\mathbb{E}\tilde\phi_{L_+}=\mathbb{E}\tilde\phi_L+\epsilon\big(\mathbb{E}\phi_L'+
\mathbb{E}\tilde\phi_{Li}\mathbb{E}\partial_i\phi_L'\big)
\quad\mbox{and}\quad\mathbb{E}\tilde\phi_{L=1}=0.
\end{align*}
Since $\mathbb{E}\phi_L'=0$, see (\ref{ao27}),
and the more obvious $\mathbb{E}\partial_i\phi_{L}'=0$ as a consequence of
the stationarity of $\phi_L'$,
we learn that centeredness is propagated, so that the first part of (\ref{ao54ter}) holds.
We now apply $\mathbb{E}$ to (\ref{ao32}) and obtain by the martingale structure,
and using the just established first part of (\ref{ao54ter})
\begin{align*}
\mathbb{E}\tilde\sigma_{L_+}=\mathbb{E}\tilde\sigma_L
+\epsilon\big(\mathbb{E}\sigma_L'+\mathbb{E}\tilde\sigma_{Li}\mathbb{E}\partial_i\phi_L'\big)\quad\mbox{and}\quad
\mathbb{E}\tilde\sigma_{L=1}=0.
\end{align*}
By $\mathbb{E}\sigma_L'=0$ and once more $\mathbb{E}\partial_i\phi_{L}'=0$ this
yields the second part of (\ref{ao54ter}).

\smallskip

We now turn to the first part of (\ref{ao64}) and apply $\mathbb{E}(\cdot)^2$ to (\ref{as31}); 
by the martingale structure, the mixed term factorizes into
$2\epsilon\mathbb{E}\phi_L'(\xi_i+\epsilon\partial_i\phi_{L}')\mathbb{E}\tilde\phi_{Li}$
and thus vanishes by (\ref{ao54ter}). The first square of (\ref{as31}) factorizes into
$\mathbb{E}(\xi_i+\epsilon\partial_i\phi_{L}')$ $(\xi_j+\epsilon\partial_j\phi_{L}')$
$\mathbb{E}\tilde\phi_{Li}\tilde\phi_{Lj}$, which according to (\ref{orthogonality}) further reduces to $\mathbb{E}|\xi+\epsilon\nabla\phi_L'|^2$
$\mathbb{E}\tilde\phi_L^2$ $=(1+\epsilon^2\mathbb{E}|\nabla\phi_L'|^2)$ 
$\mathbb{E}\tilde\phi_L^2$. Summing up, we obtain
\begin{align*}
\mathbb{E}\tilde\phi_{L_+}^2=
(1+\epsilon^2\mathbb{E}|\nabla\phi_L'|^2)\mathbb{E}\tilde\phi_L^2
+\epsilon^2\mathbb{E}{\phi_L'}^2\stackrel{(\ref{ao34bis})}{=}
(1+\frac{\epsilon^2\ln M}{2\tilde\lambda_L^2})\mathbb{E}\tilde\phi_L^2
+\epsilon^2\mathbb{E}{\phi_L'}^2.
\end{align*}
Upon multiplication by
\begin{align}\label{as34}
\tilde\lambda_{L_+}^2\stackrel{(\ref{ao30})}{=}
\big(1+\frac{\epsilon^2\ln M}{2\tilde\lambda_L^2}\big)^{2}\tilde\lambda_{L}^2,
\end{align}
this implies the iterable 
\begin{align}\label{as32}
\tilde{\lambda}_{L_{+}}^{2} \mathbb{E}\tilde\phi_{L_+}^2
=\Big(1+\frac{\epsilon^2\ln M}{2\tilde\lambda_L^2}\Big)^{3}
\tilde{\lambda}_{L}^{2} \mathbb{E}\tilde\phi_L^2
+\Big(1+\frac{\epsilon^2\ln M}{2\tilde\lambda_L^2}\Big)^{2}
\epsilon^2{\tilde\lambda_L^2}\mathbb{E}{\phi_L'}^2.
\end{align}

\smallskip

Before proceeding to the iteration, we turn to the second part of (\ref{ao64}) 
and apply $\mathbb{E}(\cdot)^2$ to (\ref{ao32}). Appealing to the 
martingale structure we obtain from the incremental centeredness $\mathbb{E}\sigma_L'$
$=\mathbb{E}\psi_L'$ $=\mathbb{E}\partial_i\phi_L'$ $=0$ that all
mixed terms between the leading $\tilde\sigma_L$ and the others vanish
besides the one coming from $\epsilon^2\tilde\phi_{Li}\psi_L'\partial_i\phi_L'$:
\begin{align}\label{ao59}
\mathbb{E}\tilde\sigma_{L_+}^2
=\mathbb{E}\tilde\sigma_{L}^2+\epsilon^2\mathbb{E}\big(
\sigma_L'+\tilde\phi_{Li}\psi_L'(\xi_i+\epsilon\partial_i\phi_L')
+\tilde\sigma_{Li}\partial_i\phi_L'\big)^2
+2\epsilon^2\mathbb{E}\tilde\phi_{Li}\tilde\sigma_L\mathbb{E}\psi_L'\partial_i\phi_L'.
\end{align}
By the cumulative centeredness (\ref{ao54ter}), two of the three
mixed terms coming from expanding the square in the middle {r.h.s.} term of (\ref{ao59}) 
vanish, and we obtain by isotropy (\ref{orthogonality})
\begin{align*}
\lefteqn{\mathbb{E}\big(\sigma_L'+\tilde\phi_{Li}\psi_L'(\xi_i+\epsilon\partial_i\phi_L')
+\tilde\sigma_{Li}\partial_i\phi_L'\big)^2}\nonumber\\
&=\mathbb{E}{\sigma_L'}^2
+\mathbb{E}\tilde\phi_{L}^2\mathbb{E}{\psi_L'}^2|\xi+\epsilon\nabla\phi_L'|^2
+\mathbb{E}\tilde\sigma_{L}^2\mathbb{E}|\nabla\phi_L'|^2
+2\mathbb{E}\tilde\phi_{Li}\tilde\sigma_{Lj}\mathbb{E}\psi_L'(\xi_i+\epsilon\partial_i\phi_L')
\partial_j\phi_L'.
\end{align*}
Since cubic terms in the Gaussian $(\phi_L',\psi_L')$ vanish
and by (\ref{as10}), this identity further simplifies to
\begin{align}\label{ao61}
\lefteqn{\mathbb{E}\big(\sigma_L'+\tilde\phi_{Li}\psi_L'(\xi_i+\epsilon\partial_i\phi_L')
+\tilde\sigma_{Li}\partial_i\phi_L'\big)^2}\nonumber\\
&=\mathbb{E}{\sigma_L'}^2
+\mathbb{E}\tilde\phi_{L}^2\big(\mathbb{E}{\psi_L'}^2
+\epsilon^2\mathbb{E}{\psi_L'}^2|\nabla\phi_L'|^2\big)
+\mathbb{E}\tilde\sigma_{L}^2\mathbb{E}|\nabla\phi_L'|^2
+2\mathbb{E}{\tilde\phi_{L}}\tilde\sigma_{Lj}\mathbb{E}\psi_L'\partial_j\phi_L'.
\end{align}
On the quartic term in (\ref{ao61}), we use Young's inequality and Gaussianity
in order to estimate $\mathbb{E}|\cdot|^4$ by $(\mathbb{E}|\cdot|^2)^2$ to the effect of
\begin{align*}
2\mathbb{E}{\psi_L'}^2|\nabla\phi_L'|^2
\le\tilde\lambda_L^{-2}\mathbb{E}{\psi_L'}^4+\tilde\lambda_L^2\mathbb{E}|\nabla\phi_L'|^4
\lesssim\tilde\lambda_L^{-2}(\mathbb{E}{\psi_L'}^2)^2
+\tilde\lambda_L^2(\mathbb{E}|\nabla\phi_L'|^2)^2
\stackrel{(\ref{ao34bis})}{\lesssim}\big(\tilde\lambda_L^{-1}\ln M\big)^2.
\end{align*}
On the non-coercive term, {i.e.} the last term in (\ref{ao61})
and likewise in (\ref{ao59}), we use Young's inequality 
and isotropy (\ref{orthogonality}):
\begin{align*}
\lefteqn{\mathbb{E}{\tilde\phi_{L}}
\tilde\sigma_{Lj}\mathbb{E}\psi_L'\partial_j\phi_L'}\nonumber\\
&\lesssim
(\tilde\lambda_L\mathbb{E}\tilde\phi_{L}^2+\tilde\lambda_L^{-1}\mathbb{E}\tilde\sigma_{L}^2)
(\tilde\lambda_L^{-1}\mathbb{E}{\psi_L'}^2+\tilde\lambda_L\mathbb{E}|\nabla\phi_L'|^2)
\stackrel{(\ref{ao34bis})}{\lesssim}
\big(\mathbb{E}\tilde\phi_{L}^2+\tilde\lambda_L^{-2}\mathbb{E}\tilde\sigma_{L}^2\big)\ln M.
\end{align*}
Inserting the last two estimates first into (\ref{ao61})
and then into (\ref{ao59}), we obtain
\begin{align}\label{ao81}
\mathbb{E}\tilde\sigma_{L_+}^2&\lesssim 
\Big(1+\frac{\epsilon^2\ln M}{\tilde\lambda_L^2} \Big)\mathbb{E}\tilde\sigma_L^2
+\Big(\frac{\epsilon^2\ln M}{\tilde\lambda_L^2}
+ \Big( \frac{\epsilon^2\ln M}{\tilde\lambda_L^2} \Big)^2\Big)
\tilde\lambda_L^2\mathbb{E}\tilde\phi_L^2+\epsilon^2\mathbb{E}{\sigma_L'}^2.
\end{align}

\smallskip

We now turn to the iteration.
Taking the sum of (\ref{as32}) and (\ref{ao81}), and noting that 
\begin{align}\label{as33}
\frac{\epsilon^2\ln M}{\tilde\lambda_L^2}
\stackrel{(\ref{ao79})}{\le}
\frac{\epsilon^2 \ln M}{1+\epsilon^2\ln L}\le 1,
\end{align}
we find
\begin{align*}
\tilde\lambda_{L_{+}}^2\mathbb{E}\tilde\phi_{L_+}^2+\mathbb{E}\tilde\sigma_{L_+}^2
\lesssim\tilde\lambda_L^2\mathbb{E}\tilde\phi_L^2+\mathbb{E}\tilde\sigma_L^2
+\epsilon^2(\tilde\lambda_{L}^2\mathbb{E}{\phi_L'}^2+\mathbb{E}{\sigma_L'}^2).
\end{align*}
Inserting (\ref{ao45}) into this,
dividing by $\epsilon^2 L_+^2$ $=M^2\epsilon^2 L^2$, cf.~(\ref{ao76}),
we obtain the iterable estimate
\begin{align*}
\frac{\tilde\lambda_{L_{+}}^2\mathbb{E}\tilde\phi_{L_+}^2+\mathbb{E}\tilde\sigma_{L_+}^2}
{\epsilon^2 L_+^2}\lesssim\frac{1}{M^2}
\frac{\tilde\lambda_L^2\mathbb{E}\tilde\phi_L^2+\mathbb{E}\tilde\sigma_L^2}{\epsilon^2 L^2}+1.
\end{align*}
For $M\gg 1$, the prefactor is $\le\frac{1}{2}$ so that iteration yields (\ref{ao64}).

\smallskip

We finally turn to \eqref{as17bis}
and apply $\mathbb{E}(\cdot)^4$ to \eqref{ao31}. Since by the martingale structure
and Gaussianity, terms involving odd powers of
$(1+\tilde\phi_{Li}\partial_i)\phi_L'$ vanish, we have
\begin{align*}
\mathbb{E}\tilde{\phi}_{L_+}^4
&=\mathbb{E}\tilde\phi_L^4
+6\epsilon^2\mathbb{E}\tilde\phi_L^2\big((1+\tilde\phi_{Li}\partial_i)\phi_L'\big)^2 
+\epsilon^4 \mathbb{E}\big((1+\tilde\phi_{Li}\partial_i)\phi_L'\big)^4,
\end{align*}
and thus obtain by Young's inequality
\begin{align*}
\mathbb{E}\tilde{\phi}_{L_+}^4
\lesssim\mathbb{E}\tilde\phi_L^4
+\epsilon^4 \mathbb{E}\big((1+\tilde\phi_{Li}\partial_i)\phi_L'\big)^4.
\end{align*}
Likewise, we have $\mathbb{E}\big((1+\tilde\phi_{Li}\partial_i)\phi_L'\big)^4$
$=\mathbb{E}{\phi_L'}^4$ $+6\mathbb{E}{\phi_L'}^2(\tilde\phi_{Li}\partial_i\phi_L')^2$
$+\mathbb{E}(\tilde\phi_{Li}\partial_i\phi_L')^4$, and thus obtain by isotropy (\ref{as11})
in form of $\mathbb{E}{\phi_{Li}'}^4$ $=\mathbb{E}{\phi_L'}^4$
\begin{align*}
\mathbb{E}\tilde{\phi}_{L_+}^4
\lesssim\mathbb{E}\tilde\phi_L^4
+\epsilon^4\big(\mathbb{E}{\phi_L'}^4+\mathbb{E}\tilde\phi_{L}^4\mathbb{E}|\nabla\phi_L'|^4\big).
\end{align*}
Appealing to Gaussianity of $\phi_L'$, we may use (\ref{ao34bis}) and (\ref{ao45}),
resulting in
\begin{align*}
\mathbb{E}\tilde{\phi}_{L_+}^4
\lesssim\Big(1+ \Big(\frac{\epsilon^2\ln M}{\tilde\lambda_L^2} \Big)^2\Big)\mathbb{E}\tilde\phi_L^4
+ \Big(\frac{\epsilon^2L_+^2}{\tilde\lambda_L^2} \Big)^2.
\end{align*}
Multiplying by $\lambda_{L_+}^4$ and dividing by $\epsilon^4 L_+^4$, 
then appealing to (\ref{as34}) and (\ref{as33}), we obtain the iterable
\begin{align*}
\frac{\tilde\lambda_{L_+}^4\mathbb{E}\tilde{\phi}_{L_+}^4}{\epsilon^4 L_+^4}
\lesssim\frac{1}{M^4}\frac{\tilde\lambda_L^4\mathbb{E}\tilde\phi_L^4}{\epsilon^4 L^4}
+1.
\end{align*}

\qed


\section{Estimate of the approximation error $f_L$} \label{s:error}

In order to pass from Lemma \ref{l:proxy} to Theorem \ref{t:main}, we need to 
show $f_L=O(\tilde\lambda_L^\frac{1}{2})$:

\begin{lemma}\label{l:estf}
For fixed but arbitrary orthogonal $R$, let $f_L^*$ denote the outcome of 
(\ref{ao18bis}) with $\xi$ replaced by $R^*\xi$. Then
\begin{align}\label{as14}
f_L^*=_{\mbox{law}}R^*f_L(R\cdot).
\end{align}
As a consequence, we have
\begin{align}\label{as35}
\mathbb{E}f_{Li}\cdot f_{Lj}
=\delta_{ij}\mathbb{E}| f_{L}|^{2}.
\end{align}
Provided $M\gg 1$ we have
\begin{align}\label{ao90}
\mathbb{E}|f_L|^2\lesssim(\epsilon^2\ln M)(1+\epsilon^2\ln M)\tilde\lambda_L.
\end{align}
\end{lemma}

{\sc Proof of Lemma \ref{l:estf}}. We start with (\ref{as14}); 
since according to (\ref{ao25}) and (\ref{ao18bis}) we have
$f_L=(1-\mathbb{E})((1+\psi_LJ)(\xi+\nabla\tilde\phi_L)-J\nabla\tilde\sigma_L)$, 
it is sufficient to show
\begin{align}\label{as16}
\big((1+\psi_L^*J)(R^*\xi+\nabla\tilde\phi_L^*),J\nabla\tilde\sigma_L^*\big)
=_{\mbox{law}}
R^*\big((1+\psi_LJ)(\xi+\nabla\tilde\phi_L),J\nabla\tilde\sigma_L\big)(R\cdot),
\end{align}
which follows from (\ref{as11}), appealing to (\ref{as15}).
Identity (\ref{as35}) follows (\ref{as14}) like (\ref{orthogonality}) followed
from (\ref{as11}).

\smallskip

We now turn to (\ref{ao90})
and apply $\mathbb{E}|\cdot|^2$ to the inductive characterization (\ref{ao55}),
where we spell out $a_L$ according to (\ref{ao25}).
By the martingale structure, we may appeal to the centeredness of all incremental
factors appearing in terms besides the second-order (in $\epsilon$) term 
$-\epsilon^2\tilde\phi_{Li}\partial_i\phi_L'J\nabla\psi_L'$, to the effect that the mixed terms vanish between the leading-order term $f_L$ and all other terms:
\begin{align*}
\lefteqn{\mathbb{E}|f_{L_+}|^2=\mathbb{E}|f_L|^2
{-2\epsilon^2\mathbb{E}\tilde\phi_{Li}f_L\cdot J
\mathbb{E}\partial_i\phi_L'\nabla\psi_L'}
+\epsilon^2\mathbb{E}\Big|\epsilon(1-\mathbb{E})\psi_L'J\nabla\phi_L'}\nonumber\\
&-\tilde\phi_{Li}(\xi_i+\epsilon\partial_i\phi_L')J\nabla\psi_L'
+\big(\tilde\phi_{Li}(1+{\epsilon\psi_LJ})-\tilde\sigma_{Li}J\big)
\nabla\partial_i\phi_L'+\partial_i\phi_L'f_{Li}\Big|^2.
\end{align*}
By the centeredness of the cumulative factors stated in  (\ref{ao54ter}), 
all the mixed terms between the purely incremental
$\epsilon(1-\mathbb{E})\psi_L'J\nabla\phi_L'$ and all the other terms vanish
besides the second-order term $\tilde\phi_{Li}\psi_LJ\nabla\partial_i\phi_L'$.
The contribution of the latter however also vanishes, 
because it is odd in the Gaussian $(\phi_L',\psi_L')$. Hence we obtain  
\begin{align}\label{ao65bis}
\lefteqn{\mathbb{E}|f_{L_+}|^2=\mathbb{E}|f_L|^2
{-2\epsilon^2\mathbb{E}\tilde\phi_{Li}f_L\cdot J\mathbb{E}
\partial_i\phi_L'\nabla\psi_L'}
+\epsilon^4\mathbb{E}|(1-\mathbb{E})\psi_L'J\nabla\phi_L'|^2}\nonumber\\
&+\epsilon^2\mathbb{E}\big|-\tilde\phi_{Li}(\xi_i+\epsilon\partial_i\phi_L')J\nabla\psi_L'
+(\tilde\phi_{Li}(1+{\epsilon\psi_{L}J})
-\tilde\sigma_{Li}J)\nabla\partial_i\phi_L'+
\partial_i\phi_L'f_{Li}\big|^2.
\end{align}
We now argue that by isotropy, also the mixed term in (\ref{ao65bis}) (and another
term of similar nature we need below) vanishes because
\begin{align}\label{ao89}
\mathbb{E}\partial_i\phi_L'\nabla\psi_L'=\mathbb{E}\partial_i\phi_L'\nabla\partial_j\phi_L'=0.
\end{align}
Indeed, for the inversion $Rx=-x$, we have that $\phi_L'^*$ considered in (\ref{as12}) 
is given by $-\phi_L'$, so that (\ref{as12}) takes the form of
\begin{align*}
(\phi_L',\psi_L')=_{\mbox{law}}(-\phi_L',\psi_L')(R\cdot),
\end{align*}
which implies
\begin{align*}
(\nabla\phi_L',\nabla^2\phi_L',\nabla\psi_L')(0)
=_{\mbox{law}}(\nabla\phi_L',-\nabla^2\phi_L',-\nabla\psi_L')(0).
\end{align*}

\smallskip

Identities (\ref{ao89}) allow us to isolate $f_L$ further:
Expanding the square of the last term in (\ref{ao65bis}), we note that
the mixed terms of $\partial_i\phi_L'f_{Li}$ with $\tilde\phi_{Lj}\xi_jJ\nabla\psi_L'$
$=\tilde\phi_LJ\nabla\psi_L'$ and with 
$(\tilde\phi_{Lj}(1+\epsilon\psi_LJ)-\tilde\sigma_{Lj}J)\nabla\partial_j\phi_L'$ vanish
by (\ref{ao89}). Also the mixed term with $\tilde\phi_{Lj}\partial_j\phi_L'J\nabla\psi_L'$ vanishes because it is cubic in the Gaussian $(\phi_L',\psi_L')$. We now turn to
the quadratic term 
\begin{align*}
\mathbb{E}\partial_i\phi_L'\partial_j\phi_L'f_{Li}\cdot f_{Lj}
=\mathbb{E}\partial_i\phi_L'\partial_j\phi_L'\mathbb{E}f_{Li}\cdot f_{Lj}
{\stackrel{(\ref{as35})}{=}\mathbb{E}|\nabla\phi_L'|^2\mathbb{E}|f_{L}|^2
\stackrel{(\ref{ao34bis})}{=}\frac{\ln M}{2\tilde\lambda_L^2}\mathbb{E}|f_{L}|^2}.
\end{align*}
Hence (\ref{ao65bis}) reduces to its final form
\begin{align}\label{ao65}
\lefteqn{\mathbb{E}|f_{L_+}|^2
=(1+\frac{\epsilon^2\ln M}{2\tilde\lambda_L^2})\mathbb{E}|f_L|^2
+\epsilon^4\mathbb{E}|(1-\mathbb{E})\psi_L'J\nabla\phi_L'|^2}\nonumber\\
&+\epsilon^2\mathbb{E}\big|-\tilde\phi_{Li}(\xi_i+\epsilon\partial_i\phi_L')J\nabla\psi_L'
+(\tilde\phi_{Li}(1+\epsilon\psi_{L}J)
-\tilde\sigma_{Li}J)\nabla\partial_i\phi_L'\big|^2.
\end{align}

\smallskip

Appealing to H\"older's inequality in probability, to the martingale structure,
and to Gaussianity in order to estimate $\mathbb{E}^\frac{1}{4}|\cdot|^4$ 
by $\mathbb{E}^\frac{1}{2}|\cdot|^2$ we have for the second {r.h.s.} term in (\ref{ao65})
\begin{align*}
\mathbb{E}|(1-\mathbb{E})\psi_L'J\nabla\phi_L'|^2
\le\big(\mathbb{E}{\psi_L'}^4\mathbb{E}|\nabla\phi_L'|^4\big)^\frac{1}{2}
\lesssim\mathbb{E}{\psi_L'}^2\mathbb{E}|\nabla\phi_L'|^2,
\end{align*}
which by (\ref{ao34bis}) yields
\begin{align}\label{ao85}
\mathbb{E}|(1-\mathbb{E})\psi_L'J\nabla\phi_L'|^2\lesssim\tilde\lambda_L^{-2}\ln^2 M,
\end{align}
and constitutes the dominant term. The last term in (\ref{ao65})
is estimated by the triangle inequality with respect to $\mathbb{E}^\frac{1}{2}|\cdot|^2$,
using the martingale structure, the isotropy (\ref{orthogonality}),
and by H\"older's and Jensen's inequality in probability, yielding
\begin{align*}
\lefteqn{\mathbb{E}^\frac{1}{2}\big|-\tilde\phi_{Li}(\xi_i+\epsilon\partial_i\phi_L')J\nabla\psi_L'
+(\tilde\phi_{Li}(1+\epsilon\psi_{L}J)
-\tilde\sigma_{Li}J)\nabla\partial_i\phi_L'\big|^2}\nonumber\\
&\le\mathbb{E}^\frac{1}{2}\tilde\phi_{L}^2\big(1+\epsilon
\mathbb{E}^\frac{1}{4}|\nabla\phi_L'|^4\big)\mathbb{E}^\frac{1}{4}|\nabla\psi_L'|^4
+\big(\mathbb{E}^\frac{1}{4}\tilde\phi_{L}^4(1+\epsilon\mathbb{E}^\frac{1}{4}\psi_{L}^4)
+\mathbb{E}\tilde\sigma_{L}^2\big)\mathbb{E}^\frac{1}{2}|\nabla^2\phi_L'|^2.
\end{align*}
Once more by Gaussianity, this improves to (after squaring the estimate)
\begin{align*}
\mathbb{E}&\big|-\tilde\phi_{Li}(\xi_i+\epsilon\partial_i\phi_L')J\nabla\psi_L'
+(\tilde\phi_{Li}(1+\epsilon\psi_{L}J)
-\tilde\sigma_{Li}J)\nabla\partial_i\phi_L'\big|^2\nonumber\\
&\lesssim\mathbb{E}\tilde\phi_{L}^2\big(1+\epsilon^2
\mathbb{E}|\nabla\phi_L'|^2\big)\mathbb{E}|\nabla\psi_L'|^2
+\big(\mathbb{E}^\frac{1}{2}\tilde\phi_{L}^4(1+\epsilon^2\mathbb{E}\psi_{L}^2)
+\mathbb{E}\tilde\sigma_{L}^2\big)\mathbb{E}|\nabla^2\phi_L'|^2,
\end{align*}
which by (\ref{ao37}), (\ref{ao34bis}), (\ref{ao54}), (\ref{ao64}), and (\ref{as17bis}) 
yields
\begin{align} \label{ao68}
\mathbb{E}&\big|-\tilde\phi_{Li}(\xi_i+\epsilon\partial_i\phi_L')J\nabla\psi_L'
+(\tilde\phi_{Li}(1+\epsilon\psi_{L}J)
-\tilde\sigma_{Li}J)\nabla\partial_i\phi_L'\big|^2\nonumber\\
&\lesssim \frac{\epsilon^2}{\tilde\lambda_L^2}
\Big(1+\frac{\epsilon^2\ln M}{\tilde\lambda_L^{2}}
+\frac{1+\epsilon^2\ln {L}}{\tilde\lambda_L^2}\Big)
{\lesssim \frac{\epsilon^2}{\tilde\lambda_L^2}
\Big(1+\frac{1+\epsilon^2\ln {L}}{\tilde\lambda_L^2}\Big)
\stackrel{(\ref{ao76})}{\lesssim}\frac{\epsilon^2}{\tilde\lambda_L^2}.}
\end{align}
Inserting the estimates (\ref{ao85}) and (\ref{ao68}) into 
the identity (\ref{ao65}) yields because of $M\gg 1$
\begin{align*}
\mathbb{E}|f_{L_+}|^2-(1+\frac{\epsilon^2\ln M}{2\tilde\lambda_L^2})\mathbb{E}|f_L|^2
\lesssim{\frac{\epsilon^4\ln^2 M}{\tilde\lambda_L^{2}}},
\end{align*}
which by (\ref{ao79}) entails 
\begin{align}\label{ao93}
\mathbb{E}|f_{L_+}|^2-\exp\big(\frac{\epsilon^2\ln M}{2(1+\epsilon^2\ln L)}\big)
\mathbb{E}|f_L|^2\lesssim\frac{\epsilon^4\ln^2 M}{1+\epsilon^2\ln L}.
\end{align}

\smallskip

By the geometric partitioning (\ref{ao76}) we obtain by iterating (\ref{ao93}) 
in $\frac{\ln L}{\ln M}$ steps
\begin{align*}
\mathbb{E}|f_{L}|^2&\lesssim(\epsilon^2\ln M)\sum_{n=0}^{\frac{\ln L}{\ln M}-1}
\frac{\epsilon^2\ln M}{1+n\epsilon^2\ln M}
\exp\big(\frac{1}{2}\sum_{m=n+1}^{\frac{\ln L}{\ln M}-1}
\frac{\epsilon^2\ln M}{1+m\epsilon^2\ln M}\big).
\end{align*}
The inner sum is estimated as in (\ref{as36}) and (\ref{as37}):
\begin{align*}
\sum_{m=n+1}^{\frac{\ln L}{\ln M}-1}
\frac{\epsilon^2\ln M}{1+ m\epsilon^2\ln M}
\le\int_n^{\frac{\ln L}{\ln M}}\dif m\frac{\epsilon^2\ln M}{1+ m\epsilon^2\ln M}
=\ln\frac{1+ \epsilon^2\ln L}{1+n\epsilon^2\ln M},
\end{align*}
so that we obtain
\begin{align*}
\mathbb{E}|f_{L}|^2& 
\lesssim(\epsilon^2\ln M)\sum_{n=0}^{\frac{\ln L}{\ln M}-1}
\frac{\epsilon^2\ln M}{1+n\epsilon^2\ln M}
\Big(\frac{1+\epsilon^2\ln L}{1+\epsilon^2n\ln M}\Big)^\frac{1}{2}.
\end{align*}
Also this (outer) sum is estimated by comparing it to an integral:
\begin{align*}
\sum_{m=0}^{\frac{\ln L}{\ln M}-1}\frac{\epsilon^2\ln M}{(1+m\epsilon^2\ln M)^\frac{3}{2}}
&\le \epsilon^2\ln M+
\int_{0}^{\infty} \dif m\frac{\epsilon^2\ln M}{(1+m\epsilon^2\ln M)^\frac{3}{2}}\nonumber\\
&=\epsilon^2\ln M+2.
\end{align*}
Hence we end up with
\begin{align*}
\mathbb{E}|f_L|^2\lesssim(\epsilon^2\ln M)(1+\epsilon^2\ln M)(1+\epsilon^2\ln L)^\frac{1}{2},
\end{align*}
so that (\ref{ao90}) follows from (\ref{ao79}). \qed


\section{Proof of Theorem \ref{t:main}}\label{s:proofmain}
In order to reach the full range $\epsilon^2<\infty$ for the
amplitude $\epsilon$ of the divergence-free drift, we use that the effective
diffusivity is non-decreasing in the squared amplitude
-- a classical observation, which we reproduce for the reader's convenience:

\begin{lemma}\label{l:monotone}
$\lambda_{L}$ is non-decreasing in\footnote{Since $-\psi$ and $\psi$ have the same law,
$\lambda_L$ is an even function in $\epsilon$, and thus a function of $\epsilon^2$ only.} 
$\epsilon^2$.
\end{lemma}

{\sc Proof of Lemma \ref{l:monotone}}. 
Formally, by (\ref{ao25}), the $\epsilon$-derivative of (\ref{ao01}) is given by
\begin{align*}
-\nabla\cdot\big(a_L\frac{\partial\nabla\phi_{L}}{\partial\epsilon}
+\psi_L J(\xi+\nabla \phi_{L})\big)=0
\quad\mbox{and}\quad\mathbb{E}\frac{\partial\nabla\phi_L}{\partial\epsilon}=0.
\end{align*}
Once more by (\ref{ao25}) and (\ref{ao01}), and for $\epsilon\not=0$, 
this PDE can be rewritten as
\begin{align}\label{derphieps}
-\nabla\cdot\big(a_L\frac{\partial\nabla\phi_{L}}{\partial\epsilon}
-\frac{1}{\epsilon}\nabla \phi_{L}\big)=0,
\end{align}
which rigorously defines the derivative $\frac{\partial\nabla\phi_{L}}{\partial\epsilon}$
as a square-integrable, centered and curl-free vector field. The latter is 
an admissible test function for (\ref{derphieps}), yielding in view of (\ref{ao40})
\begin{equation*}
\mathbb{E}\big|\frac{\partial\nabla\phi_{L}}{\partial\epsilon}\big|^2  
=\frac{1}{\epsilon}\mathbb{E}\nabla \phi_{L}
\cdot\frac{\partial\nabla\phi_{L}}{\partial\epsilon}.
\end{equation*}
We then learn from (\ref{ao08}) that $\lambda_L$ is differentiable in $\epsilon\not=0$
with
\begin{align*}
\frac{1}{2\epsilon}\frac{\partial\lambda_L}{\partial\epsilon}
\stackrel{(\ref{ao08})}{=}\frac{1}{\epsilon}\mathbb{E}\nabla \phi_{L}
\cdot\frac{\partial\nabla\phi_{L}}{\partial\epsilon}
=
\mathbb{E}\big|\frac{\partial\nabla\phi_{L}}{\partial\epsilon}\big|^2\ge 0.
\end{align*}
\qed

{\sc Proof of Theorem \ref{t:main}}. As long as $L\sim 1$, (\ref{ao38}) is trivial,
and the asymptotics $C_\epsilon\sim 1$ for $\epsilon\ll 1$ amount to a classical perturbation
result. Hence without loss of generality we assume $L\gg 1$; 
then there exists an $M$ with $\frac{\ln L}{\ln M}\in\mathbb{N}$ 
so that we are (exactly) in the geometric setting of 
(\ref{ao76}){, and that is sufficiently large such that (\ref{ao90}) holds. 
On the one hand, by (\ref{ao79}) in Lemma \ref{l:proxy}  (and using $M\sim 1$)
we have for the proxy $\tilde\lambda_L$
\begin{align*}
|\tilde\lambda_L^2-(1+\epsilon^2\ln L)|\lesssim{\epsilon^4\ln L}
{\quad\mbox{and thus}\quad
\big|\frac{\tilde\lambda_L^2}{1+\epsilon^2\ln L}-1\big|\lesssim\epsilon^2}.
\end{align*}
On the other hand by (\ref{ao80}) in the same lemma and by (\ref{ao90}) in Lemma \ref{l:estf}
(and using $M\sim 1$) we obtain
\begin{align*}
|\lambda_L-\tilde\lambda_L|\lesssim|\epsilon|\tilde\lambda_L
\quad\mbox{and thus}\quad
{\big|\frac{\lambda_L^2}{\tilde\lambda_L^2}-1\big|\lesssim}{|\epsilon|}.
\end{align*}

This establishes (\ref{ao38}) for $\epsilon^2\ll 1$, including the asymptotics, 
but only the upper bound for $\epsilon^2\gtrsim 1$. The missing lower bound for $\epsilon^2\gtrsim 1$
then follows from the monotonicity of Lemma \ref{l:monotone}.

\qed 

\section{Proof of Theorem \ref{t:laplace}} \label{s:proofmain2}
We first state the following representation lemma, which allows us to express $\mathbb{E}|X_t|^2$ as the energy of the solution of a parabolic equation. We believe this observation is already present in literature, but we are unable to find it. Therefore we present both the statement and proof for the reader's convenience.
	\begin{lemma}\label{lem:repfmla}
		Let $X_t$ satisfy the SDE \eqref{e:sde}, and $v=v(t,x)$ be the stationary solution of the parabolic equation
		\begin{equation}\label{rr08}
			\partial_t v - \Delta v-\epsilon b\cdot (\xi+\nabla v) = 0, \, v(t=0,\cdot)=0,
		\end{equation}
	then for any unit vector $\xi$,
	\begin{equation}\label{rr01}
		\frac{1}{2}\mathbb{E}(\xi\cdot X_T)^2 = T + \mathbb{E} \Big(v(T)^2 + \int_0^T \dif t |\nabla v|^2\Big).
	\end{equation}
	\end{lemma}

\smallskip

Note that a standard truncation argument combined with energy-type estimates implies that there exists a unique stationary Lax-Milgram solution $v$ to \eqref{rr08} and that both $v$ and $\nabla v$ are bounded for any algebraic moment.

\smallskip
 
{\sc Proof of Lemma \ref{lem:repfmla}.}
	The {l.h.s.} of (\ref{rr01}) can be rephrased in terms
	of the quenched heat kernel $G(T,x,y)$ as
	\begin{align}\label{rr02}
	\frac{1}{2}\mathbb{E}(\xi\cdot X_T)^2
		=\frac{1}{2}\int \dif x(\xi\cdot x)^2 \mathbb{E} G(T,x,0).
	\end{align}
	In view of its characterization
	\begin{align}\label{rr03}
		\partial_tG-\Delta_xG-\epsilon b(x) \cdot \nabla_x G=0\quad\mbox{and}\quad
		G(t=0,x,y)=\delta(x-y)
	\end{align}
	we obtain from an integration by parts in $x$ in (\ref{rr02}) and the fact that $b$ is divergence free
	\begin{align*}
\frac{1}{2}	\mathbb{E}(\xi\cdot X_T)^2
		&=\frac{1}{2}\mathbb{E}\int_0^T \dif t \int \dif x(\xi\cdot x)^2
		\partial_tG(t,x,0)\nonumber\\
		&=\mathbb{E}\int_0^T\dif t \int \dif x
		\big(1-(\xi\cdot x)(\epsilon\xi\cdot b(x))\big)G(t,x,0).
	\end{align*}
	Since (\ref{rr03}) implies in particular $\int \dif x G(t,x,y)=1$,
	in order to establish (\ref{rr01}), it remains to show
	\begin{align}\label{rr04}
		-\mathbb{E}\int_0^T \dif t\int \dif x
		(\xi\cdot x)(\epsilon\xi\cdot b(x))G(t,x,0)= \mathbb{E} \Big(v(T)^2 + \int_0^T \dif t |\nabla v|^2\Big).
	\end{align}
	Multiplying (\ref{rr08}) with $v$, applying $\mathbb{E}$, integrating in time and appealing to
	stationarity of $v(T)$ and $b$, we obtain from the fact that the latter is
	divergence-free
	\begin{align*}
	\mathbb{E} \Big(v(T)^2 + \int_0^T\dif t |\nabla v|^2\Big)
		=\mathbb{E}(\epsilon\xi\cdot b)\int_0^T \dif t \, v,
	\end{align*}
	so that removing the $\int_0^T \dif t$ integral, (\ref{rr04}) can be rephrased as
	\begin{align*}
		-\mathbb{E}\int \dif x
	(\xi\cdot x)(\epsilon \xi\cdot b(x))G(T,x,0)=\mathbb{E}(\epsilon\xi\cdot b)v(T).
	\end{align*}
	By the pathwise representation of $v(T)$, cf.~(\ref{rr08}), in terms  
	of the parabolic Green function, this can be rewritten as 
	\begin{align}\label{rr05}
		-\mathbb{E}\int \dif x
		(\xi\cdot x)(\epsilon\xi\cdot b(x))G(T,x,0)=\mathbb{E}(\epsilon\xi\cdot b(0))\int_0^T \dif t  \int \dif y(\epsilon\xi\cdot b(y))G(t,0,y).
	\end{align}
	
	\smallskip
	
	In order to establish (\ref{rr05}), we use that the quenched parabolic Green function 
	$G(b,t,x,y)$
	satisfies the dual equation in the $y$-variable, and that by uniqueness it is shift invariant:
	\begin{align}
		\partial_t G(t, x,\cdot)-\Delta_yG(t, x,\cdot)+\epsilon \nabla \cdot G(t,x,\cdot) b(y) = 0,
		\label{rr12}\\
		G(b(\cdot+z),t,x,y)=G(b,t,x+z,y+z).\label{rr11}
	\end{align}
	This allows us to conclude: By stationarity of $\mathbb{E}$, the {l.h.s.} of (\ref{rr05}) 
	can be rewritten as $	\mathbb{E}\int \dif x
	(\xi\cdot x)(\epsilon \xi\cdot b(x+z))G(b(\cdot+z),T,x,0)$,
	where $\mathbb{E}$ acts on $b$. Using (\ref{rr11}) for $z=-x$ (and with $y=0$), we see that
	this is identical to $	\mathbb{E}\int \dif x
	(\xi\cdot x)(\epsilon \xi\cdot b(0))G(T,0,-x)$,
	which we write as $-\mathbb{E}(\epsilon\xi\cdot b(0)) \int \dif y(\xi\cdot y)G(T,0,y)$.
	However, we obtain from testing (\ref{rr12}) (with $x=0$) with the linear function $\xi\cdot y$,
	using integration by parts,
	\begin{align*}
		\int \dif y(\xi\cdot y)\partial_t G(t,0,y)=\int \dif y(\epsilon\xi\cdot b(y))G(t,0,y).
	\end{align*}
	Integrating in time, this establishes (\ref{rr05}) and thus (\ref{rr01}). \qed

\smallskip

We now estimate the energy of $v(T)$ shown in \eqref{rr01}. Our strategy is to compare $v(T)$ with $\tilde{\phi}_L$ for $T=L^2$ and use energy estimates. The reason we choose $T=L^2$ is the parabolic scaling.

\smallskip

{\sc Proof of Theorem \ref{t:laplace}}. Let us fix $L=\sqrt{T} \geq 1$. We first apply $\nabla\cdot$ to (\ref{ao18bis}) and take the difference with \eqref{rr08}. Then using that $\tilde\phi_L$ is time-independent, we finally get that $v-\tilde{\phi}_L$ solves
\begin{equation}\label{vminphiL}
		\left\{\begin{aligned}& \partial_t (v-\tilde{\phi}_L) -\Delta (v-\tilde{\phi}_L) -\epsilon b\cdot \nabla (v-\tilde{\phi}_L) = \epsilon (b-b_L)\cdot \xi - \epsilon \nabla\cdot \tilde{\phi}_L (b-b_L)-\nabla\cdot f_L, \\ & (v-\tilde{\phi}_L)(t=0)= -\tilde{\phi}_L. \end{aligned}\right.
	\end{equation}
Next we test \eqref{vminphiL} against $v-\tilde{\phi}_L$ which, using that $\tilde{\phi}_L$ is independent of $b-b_L$, yields
\begin{align*}
	\frac{1}{2}\frac{\dif}{\dif t} \mathbb{E}(v-\tilde{\phi}_L)^2 & +  \mathbb{E}|\nabla(v-\tilde{\phi}_L)|^2 \\  & \lesssim \epsilon^2 L^2 \mathbb{E} |b-b_L|^2 + \frac{1}{L^2} \mathbb{E}(v-\tilde{\phi}_L)^2 +\epsilon^2\mathbb{E}\tilde{\phi}_L^2 \mathbb{E}|b-b_L|^2  + \mathbb{E}|f_L|^2. 
\end{align*}
Note that we may bound the difference $\mathbb{E} |b-b_L|^2$ by invoking \eqref{e:bminbL}, then making use of the moment bounds \eqref{ao64} on the proxy $\tilde{\phi}_L$ and \eqref{ao90} on the error $f_L$ we may estimate
\begin{align*}
	\frac{1}{2}\frac{\dif}{\dif t} \mathbb{E}(v-\tilde{\phi}_L)^2  +  \mathbb{E}|\nabla(v-\tilde{\phi}_L)|^2 \lesssim  \frac{1}{L^2} \mathbb{E}(v-\tilde{\phi}_L)^2  + \epsilon^2\tilde{\lambda}_L.
\end{align*}
Hence Gronwall's inequality implies
\begin{align*}
	 \mathbb{E}(v(T)-\tilde{\phi}_L)^2 + \int_0^T \dif t \, e^{(T-t)/L^2}  \mathbb{E}|\nabla(v-\tilde{\phi}_L)|^2  & \lesssim e^\frac{T}{L^2} \mathbb{E} \tilde\phi_L^2 + \int_0^T \dif t \, e^{(T-t)/L^2} \epsilon^2 \tilde\lambda_L,
\end{align*}
which after employing \eqref{ao64} again and recalling that $T=L^2$, gives
\begin{equation}\label{engydiff} 
	 \mathbb{E}(v(T)-\tilde{\phi}_L)^2 + \int_0^T \dif t \,  \mathbb{E}|\nabla(v-\tilde{\phi}_L)|^2  \lesssim  \epsilon^2 \tilde{\lambda}_L T.
\end{equation}

\smallskip 

Let us now define the energy
	\begin{equation*}
		\mathcal{E}(T) := 1 + \frac{1}{T} \mathbb{E}\Big(v(T)^{2} + \int_{0}^{T}\dif t |\nabla v|^{2}  \Big) \stackrel{\eqref{rr01}}{=} \frac{1}{2T}\mathbb{E}(\xi\cdot X_T)^2.
	\end{equation*}
 Since $\tilde{\phi}_{L}$ is independent of time, we obtain from expanding the square
\begin{equation}\label{ETexp}
	\begin{aligned}
		\mathcal{E}(T) &= 1 + \frac{1}{T} \mathbb{E}\tilde{\phi}_{L}^{2} + \mathbb{E}|\nabla \tilde{\phi}_{L}|^{2} + \frac{2}{T} \mathbb{E}(v(T) - \tilde{\phi}_{L})  \tilde{\phi}_{L} + \frac{2}{T} \mathbb{E}\int_{0}^{T} \dif t \, \nabla (v - \tilde{\phi}_{L}) \cdot \nabla \tilde{\phi}_{L}   \\
			&\qquad + \frac{1}{T} \mathbb{E}(v(T) - \tilde{\phi}_{L})^{2} + \frac{1}{T} \mathbb{E}\int_{0}^{T}\dif t \, |\nabla (v - \tilde{\phi}_{L})|^{2} .
	\end{aligned} \end{equation}
In order to estimate the {r.h.s.} of \eqref{ETexp} we first observe that we need a bound for $\nabla\tilde\phi_L$. From the upper bound \eqref{ao90} of $f_L$, we may turn the upper bound \eqref{ao97} into
 \begin{equation}\label{newbdgrphi}
     \mathbb{E}|\nabla\tilde\phi_L|^2 \lesssim \tilde \lambda_L + \epsilon^2\tilde\lambda_L \lesssim \tilde \lambda_L.
 \end{equation}
 Recalling the identity \eqref{ao96} and combining again with \eqref{ao90}, we may control the leading order term in \eqref{ETexp} by
 \begin{equation}\label{mainET}
     \left|1+\mathbb{E}|\nabla \tilde{\phi}_{L}|^{2} - \tilde\lambda_L\right| = |\mathbb{E} f_L\cdot\nabla\tilde\phi_L| \le \mathbb{E}^\frac{1}{2}|f_L|^2\mathbb{E}^\frac{1}{2}|\nabla\tilde\phi_L|^2 \lesssim |\epsilon| \tilde\lambda_L.
 \end{equation}
Now for the cross terms in \eqref{ETexp}, we use first the Cauchy-Schwarz inequality and then combine with  \eqref{ao64}, \eqref{engydiff}, \eqref{newbdgrphi} and the fact that $\tilde\lambda_L\ge 1$ to finally get
	\begin{align*}
		\frac{2}{T}\Big|\mathbb{E}(v(T) & - \tilde{\phi}_{L}) \tilde{\phi}_{L}  +\mathbb{E}\int_{0}^{T} \dif t \,\nabla (v - \tilde{\phi}_{L}) \cdot \nabla \tilde{\phi}_{L}   \Big| \\  & \leq \frac{2}{T} \mathbb{E}^\frac{1}{2}(v(T) - \tilde{\phi}_{L})^{2}  \mathbb{E}^\frac{1}{2}\tilde{\phi}_{L}^{2}  + \frac{2}{\sqrt{T}} \mathbb{E}^\frac{1}{2}\left(\int_{0}^{T} \dif t|\nabla (v - \tilde{\phi}_{L})|^{2}\right) \mathbb{E}^\frac{1}{2}|\nabla \tilde{\phi}_{L}|^{2}  \\
			& \lesssim\sqrt{\epsilon^{2} \tilde{\lambda}_{L} \cdot \epsilon^{2} \tilde{\lambda}_{L}^{-2}} +  \sqrt{\epsilon^2\tilde\lambda_L^2} \lesssim |\epsilon| \tilde \lambda_L. \stepcounter{equation} \tag{\theequation} \label{ETcross}
	\end{align*}

Putting all the above together we are now in position to estimate $\mathcal{E}(T)$ via \eqref{ETexp}. Indeed, the first term could be estimated with the aid of \eqref{ao64}, while for the second term we may apply \eqref{mainET}. Next observe that, in view of \eqref{engydiff}, the second line of \eqref{ETexp} is bounded above by $\epsilon^{2} \tilde{\lambda}_{L}$. Combining with \eqref{ETcross} for the cross terms, we get the final estimate
 \begin{align*}
     |\mathcal{E}(T) -\tilde\lambda_L|\lesssim |\epsilon|\tilde\lambda_L.
 \end{align*}
Combined with \eqref{ao79}, this establishes  \eqref{eq:EXT2} for $\epsilon^2\ll 1$ including the asymptotics, and the upper bound for $\epsilon^2\gtrsim 1$. The lower bound for  $\epsilon^2 \gtrsim 1$ is again recovered from monotonicity of $\mathcal{E}(T)$ in $\epsilon^2$, following an identical argument to Lemma \ref{l:monotone}.

\begin{appendix}

\section{Heuristic Derivation of the Scaling Law} \label{sec:appendix_heuristic}
We add more details on the heuristics from \cite{FFQSS85} for the reader's convenience. The idea of the renormalization group argument is to perturb around
the standard diffusion. This is implemented by rescaling the spatial coordinates
in a time-dependent way such that always $\mathbb{E}|\hat X(t)|^2=t$.
Given a time instance $t$, we consider a time interval $[t,t_+)$ over
which the second moment increases by a factor of $M^2$.
Hence given $t$ and $M$, we define $t_+$ and the exponent $z$ via
\begin{align}\label{ab02}
\mathbb{E}|X_{t_+}|^2=M^2\mathbb{E}|X_t|^2\quad\mbox{and}\quad t_+=M^zt.
\end{align}
We expect $z$ to be $<2$ and  depend on $t$.
In order to propagate the assumption 
\begin{align}\label{ab03}
\mathbb{E}|X_t|^2=t,
\end{align}
to the next time instance $t_+$ in form of
$\mathbb{E}|\hat X_{t_+}|^2=t_+$, we need to (incrementally) rescale space,
and thus the drift (which has units of velocity), and its covariance as (which
is quadratic in the drift) as
\begin{align}\label{ab04}
x=M^{1-\frac{z}{2}}\hat x\quad\mbox{and}\quad b=M^{1-\frac{z}{2}}\hat b
\quad\mbox{so that}\quad \hat R=M^{z-2} R.
\end{align}
In particular, we expect the effect of the drift to become smaller
under this incremental rescaling.

\smallskip

Now comes the core of the argument: We think of $0<M-1\ll 1$, so that
we may expect $\frac{\dif}{\dif s}\mathbb{E}|X_s|^2$ to be approximately constant
in the (logarithmically small) time interval $[t,t_+)$. Because the drift
is divergence-free (its generator is skew symmetric in the physics jargon), 
it increases the effective diffusivity. Since the drift is independent from
the Brownian motion, the effect is additive. Moreover, by symmetry, the effect must
be to leading order quadratic in $b$, or linear in its covariance $R$. 
Neglecting the UV cut-off, the covariance function $R$
is universal, namely a multiple $T>0$ of the kernel of the Leray projection, that is
\begin{align}\label{ab09}
R(x)=T J\frac{1}{|x|^2}({\rm id}-2\frac{x}{|x|}\otimes\frac{x}{|x|})J^t.
\end{align}
Hence there is a universal constant $c_0>0$ such that
\begin{align*}
\frac{\dif}{\dif s}\mathbb{E}|X_s|^2\approx 1+c_0 T.
\end{align*}
Integrating it over $[t,t_+)$ and using (\ref{ab02}) and (\ref{ab03}),
this yields
\begin{align*}
(M^2-1)t\approx (1+c_0T)(M^z-1)t.
\end{align*}
Because of $0<M-1\ll 1$, we learn
\begin{align}\label{ab05}
2\approx z(1+c_0 T)\quad\mbox{and thus}\quad c_0 T\approx\frac{2}{z}-1.
\end{align}
In terms of the form (\ref{ab09}), the rescaling (\ref{ab04}) translates into
\begin{align*}
\hat T=M^{2(z-2)} T,
\end{align*}
which yields
\begin{align*}
\frac{\hat T-T}{M-1}\approx2(z-2)T\stackrel{(\ref{ab05})}{\approx} -2\frac{(2-z)^2}{z}.
\end{align*}
Since the relation (\ref{ab05}) also holds for the exponent for
the next increment, that is $\frac{2}{\hat z}-1\approx \hat T$, we obtain
a discrete dynamical system for $z$:
\begin{align}\label{ab06}
\frac{1}{M-1}\big(\frac{1}{\hat z}-\frac{1}{z}\big)\approx -\frac{(2-z)^2}{z}.
\end{align}

\smallskip

It remains to leverage (\ref{ab06}) by calculus methods.
In view of the second item in (\ref{ab02}) in form of $\ln t_+-\ln t$ $=z\ln M$
$\approx z(M-1)$, this discrete dynamical system can be assimilated to a continuum one
in logarithmic time:
\begin{align*}
z\frac{\dif}{\dif\ln t}\frac{1}{z}\approx-\frac{(2-z)^2}{z},\quad\mbox{that is}\quad
\frac{\dif z}{\dif \ln t}\approx(2-z)^2.
\end{align*}
All solutions $z<2$ of this autonomous ODE are parameterized by 
some characteristic time scale $\tau$ via
\begin{align}\label{ab07}
z\approx2-\frac{1}{\ln\frac{t}{\tau}}\quad\mbox{for}\;t\gg\tau.
\end{align}
On the other hand, the continuum version of (\ref{ab02}) reads
\begin{align*}
\frac{\dif\ln\mathbb{E}|X_t|^2}{\dif\ln t}=\frac{2}{z}.
\end{align*}
Combining this with (\ref{ab07}) we obtain
\begin{align*}
\frac{\dif\ln\mathbb{E}|X_t|^2}{\dif \ln\frac{t}{\tau}}\approx\frac{2}{2-\frac{1}{\ln\frac{t}{\tau}}}
\quad\mbox{for}\;t\gg\tau,
\end{align*}
and thus by plain integration for some characteristic length scale $\ell$
\begin{align*}
\ln\frac{\mathbb{E}|X_t|^2}{\ell^2}\approx\ln\frac{t}{\tau}
+\frac{1}{2}\ln(\ln\frac{t}{\tau}-\frac{1}{2})
\quad\mbox{for}\;t\gg\tau,
\end{align*}
which by exponentiating implies the desired scaling law:
\begin{align*}
\mathbb{E}|X_t|^2\approx\ell^2\frac{t}{\tau}\sqrt{\ln\frac{t}{\tau}}
\quad\mbox{for}\;t\gg\tau.
\end{align*}

\section{Preliminaries from Homogenization Theory} \label{sec:appendix_prelim}

In this appendix, we recall two relevant facts from homogenization theory, although they are technically not necessary for any of the results proved herein.  First, we recall the connection between the effective diffusivity as defined by $\lambda_{L}$ and the asymptotics of the diffusion process $X^{(L)}_{t}$.  Next, we provide the relevant existence and uniqueness result for the flux corrector $\sigma_{L}$.

\subsection{The Effective Diffusivity} In this subsection, we recall the connection between the effective diffusivity $\lambda_{L}$ and the asymptotics of the associated diffusion process $X^{(L)}_{t}$ defined in the introduction.  In particular, we sketch the proof of the identity \eqref{e:diffusivity-X_L}.  

\smallskip

First, observe that, for a given $\xi \in \mathbb{R}^{2}$, the corrector $\phi_{L}$ is defined in precisely such a way that the process $M_{t}^{(L)}$, defined below, becomes a martingale:
	\begin{align} \label{e:martingale}		
		M_{t}^{(L)} = \xi \cdot X^{(L)}_{t} + \phi_{L}(X_{t}^{(L)}). 
	\end{align}
Indeed, in order for this to be a martingale, the function $\xi \cdot x + \phi_{L}(x)$ needs to be in the kernel of the generator of $X^{(L)}_{t}$, which follows from \eqref{ao01}.

\smallskip

%
%
%
Without going into the details here, one can prove, using the so-called ``environment seen by the particle," that 
	\begin{align*}
		\lim_{t \to \infty} \frac{\mathbb{E} (M^{(L)}_{t})^{2}}{2t} = \mathbb{E} |\xi + \nabla \phi_{L}|^{2}.
	\end{align*}
At the same time, since $\phi_{L}$ is sublinear at infinity (having a stationary, mean-zero gradient), one expects that the correction term has a lower-order contribution in the limit, and thus 
	\begin{align*}
		\lim_{t \to \infty} \frac{ \mathbb{E} (\xi \cdot X_{t}^{(L)}  )^{2} }{2t} = \lim_{t \to \infty} \frac{ \mathbb{E} (M_{t}^{(L)})^{2} }{2t} = \mathbb{E} |\xi + \nabla \phi_{L\xi} |^{2}.
	\end{align*}
All of this is rigorously justified in \cite[Corollary 11.5]{KLO12}.  Combining the last equality with \eqref{ao08}, we conclude that $\lambda_{L}$ captures the asymptotic diffusion constant of $X^{(L)}_{t}$ as in \eqref{e:diffusivity-X_L}.


\subsection{Flux Corrector $\sigma_{L}$} Next, we recall that, in homogenization error estimates, it is convenient to define the so-called flux corrector $\sigma_{L}$ to go alongside the gradient corrector $\phi_{L}$.  In the context of stochastic homogenization, this was first realized in \cite{GNO20}.  The flux corrector $\sigma_{L}$ is a random function with stationary, mean-zero gradient $\nabla \sigma_{L}$ such that
	\begin{equation}
		a_{L} (\xi + \nabla \phi_{L}) = \lambda_{L} \xi + J \nabla \sigma_{L}. \label{e:flux_correction}
	\end{equation}
As is readily checked, the equation \eqref{e:flux_correction} reflects the fact that the flux $a_{L}(\xi + \nabla \phi_{L})$ is divergence free; and $\nabla \sigma_{L}$ has vanishing expectation since $\mathbb{E} a_{L} (\xi + \nabla \phi_{L}) = \lambda_{L} \xi$. 

	\begin{theorem} For any $L \geq 1$ and any $\xi \in \mathbb{R}^{2}$, there is a unique stationary, mean-zero gradient $\nabla \sigma_{L}$ such that \eqref{e:flux_correction} holds in $\mathbb{R}^{2}$. \end{theorem}

	\begin{proof} Applying the curl operator $\nabla \cdot J$ to both sides of \eqref{e:flux_correction}, one finds
		\begin{align*}
			- \Delta \sigma_{L} = \nabla \cdot J a_{L}(\xi + \nabla \phi_{L}).
		\end{align*}
	This is precisely the equation treated in \cite[Lemma 1]{BFO18} and \cite[Proposition 2.4]{Fe20}.  \end{proof}
	
%

\end{appendix}

\begin{acks}[Acknowledgments]
GC, FO and LW would like to express their special thanks to Jing An for many 
interesting conversations and idea-exchanging they had on the topic during 
her one-year stay at MPI. FO likes to thank Nicolas Perkowski for
explaining the method of \cite{CHT22} to him.  PM thanks Takis Souganidis for helpful discussions.
\end{acks}

\begin{funding}
The second author was supported by NSF Grant DMS-2202715.
\end{funding}

\bibliographystyle{imsart-number} 
\bibliography{biblio}       

%
%
%

\end{document}